\documentclass[11pt,a4paper,reqno]{amsart}
\usepackage[T1]{fontenc}
\usepackage{lmodern,amsmath,amssymb,mathtools}
\usepackage[margin=28mm]{geometry}
\usepackage{needspace}
\usepackage[hidelinks]{hyperref}

\numberwithin{equation}{section}
\newtheorem{theorem}{Theorem}[section]
\newtheorem{lemma}[theorem]{Lemma}
\newtheorem{proposition}[theorem]{Proposition}
\newtheorem{numlemma}{Numerical Lemma}[section]
\newcommand{\R}{\mathbb R}
\newcommand{\Z}{\mathbb Z}
\newcommand{\one}{\mathbf 1}

\title[The Graf--Luschgy conjecture]
{The BCC lattice is the optimal 3D lattice quantizer for all $p\geq1$}
\author{Senping Luo}
\hypersetup{pdfauthor={Senping Luo},
pdftitle={The BCC lattice is the optimal 3D lattice quantizer for all $p\geq1$}}
\subjclass[2020]{52C15, 94A29, 65G40}
\keywords{Lattice quantization, Voronoi moments, body-centered cubic lattice, Selling parameters}

\begin{document}
\begin{abstract}
For every finite $p\ge1$, the body-centered cubic lattice uniquely
minimizes the $p$th moment of the Voronoi cell among three-dimensional
lattices of covolume one, up to orthogonal transformations. This
proves the Graf--Luschgy conjecture for three-dimensional lattice
quantizers (\cite{GL}).
\end{abstract}
\maketitle

\section{Introduction}
Following Bourne, Cagnetti and Roper~\cite[Eqs.~(1.1)--(1.4)]{BCR},
let $A=(t_1\mid t_2\mid t_3)\in GL(3,\R)$ and write
\[
 \Lambda(A)=\operatorname{span}_{\Z}\{t_1,t_2,t_3\}=A\Z^3,\qquad
 V(A)=\{x\in\R^3:\ |x|\le |x-w|\ \text{for every }w\in\Lambda(A)\}.
\]
The volume $|V(A)|=|\det A|$ is the covolume of $\Lambda(A)$.
For $p\ge1$ define the normalized Voronoi moment by
\[
 E_p(A)=|V(A)|^{-1-p/3}\int_{V(A)}|x|^p\,dx .
\]
For $c>0$, $O\in O(3)$ and $U\in GL(3,\Z)$,
$V(cOAU)=cOV(A)$, so
\begin{equation}\label{eq:basis-invariance}
 E_p(cOAU)
 = (c^3|V(A)|)^{-1-p/3}c^{p+3}\int_{V(A)}|x|^p\,dx
 =E_p(A).
\end{equation}
Taking $c=|V(A)|^{-1/3}$ reduces the minimization to $|V(A)|=1$.

Use the BCC basis of~\cite[Eq.~(4.13)]{BCR},
\[
 A_{\mathrm{bcc}}=
 \begin{pmatrix}
  1&0&1/2\\
  0&1&1/2\\
  0&0&1/2
 \end{pmatrix},
 \qquad D_3^*=\Lambda(A_{\mathrm{bcc}}).
\]
Here $|V(A_{\mathrm{bcc}})|=1/2$. For a lattice $L=\Lambda(A)$,
write $|L|=|V(A)|$ and $E_p(L)=E_p(A)$. This value is independent
of the choice of basis by~\eqref{eq:basis-invariance}. Set
\[
 L_{\rm B}=2^{1/3}D_3^*=\Lambda(2^{1/3}A_{\mathrm{bcc}}),
 \qquad |L_{\rm B}|=1.
\]

\begin{theorem}\label{thm:main}
Let $1\le p<\infty$. Among full-rank lattices $L\subset\R^3$
of covolume one,
\begin{equation}\label{eq:minimizers}
 \min_{|L|=1}E_p(L)=E_p(L_{\rm B})=E_p(A_{\mathrm{bcc}}),
 \qquad
 \operatorname*{argmin}_{|L|=1}E_p(L)=\{O L_{\rm B}:O\in O(3)\}.
\end{equation}
\end{theorem}
This resolves the Graf-Luschgy Conjecture \cite[Conjecture~8.17]{GL}.
The normalization of the sharp constant is
\[
 E_p(L_{\rm B})
 =\int_{V(2^{1/3}A_{\mathrm{bcc}})} |x|^p\,dx
 =2^{1+p/3}\int_{V(A_{\mathrm{bcc}})} |x|^p\,dx .
\]
The fixed-square expression for this constant is given in~\eqref{eq:bcc-value}.

Let $X$ be a nonempty finite set in a domain $\Omega$ of
volume one. With the masses on $X$ free, the transport cost is
\[
 \min_{\substack{f:\Omega\to X\\f\ \mathrm{measurable}}}
 \int_\Omega |x-f(x)|^p\,dx
 =\int_\Omega\min_{w\in X}|x-w|^p\,dx
 =\int_\Omega\operatorname{dist}(x,X)^p\,dx .
\]
A nearest-site map attains the minimum~\cite[Lemma~2.2]{BCR}.
On the unit cube in $\R^d$, the asymptotic quantization theorem
gives~\cite{Zador,BW}
\[
 \lim_{N\to\infty}N^{p/d}\!
 \inf_{\substack{X\subset[0,1]^d\\\#X=N}}
       \int_{[0,1]^d}\min_{w\in X}|x-w|^p\,dx\in(0,\infty).
\]
Gruber extended this asymptotic theory to Riemannian manifolds
and normed spaces~\cite{GruberOpt}. For a bounded open
$\Omega\subset\R^3$ of volume one with
$\partial\Omega=\varphi(K)$, where $K\subset\R^2$ is compact
and $\varphi:\R^2\to\R^3$ is Lipschitz,
Bourne, Cagnetti and Roper~\cite[Theorem~A.1]{BCR} prove the
pointwise limit
\[
 \lim_{\varepsilon\downarrow0}
  \varepsilon^{-p}|V(A)|^{-p/3}\,
    \int_\Omega\operatorname{dist}(x,\Lambda(\varepsilon A)\cap\Omega)^p\,dx
   =E_p(A).
\]
The same theorem gives $\Gamma$-convergence of the rescaled
distortion on $GL(3,\R)$ with the Frobenius metric.

In dimension two the moment theorem of Fejes T\'oth yields
triangular-lattice optimality for every finite $p\ge1$.
The regular hexagon also gives the unrestricted planar
asymptotic constant~\cite{FT,Gruber}. See~\cite[Sect.~2]{BCR}.
In dimension three the unrestricted problem is the setting of
Gersho's conjecture. It concerns the asymptotic cell geometry of unconstrained
optimal point sets~\cite{Gersho}. For a precise quadratic
formulation see Choksi and Lu~\cite[Conjecture~1.1]{CL}.
Du and Wang obtained numerical evidence for the BCC cell in
the three-dimensional quadratic problem~\cite{DW}.
For the high-rate variable-rate problem on periodic tilings,
Sloane and Vaishampayan computed the asymptotic constant in
terms of the cell volumes and second moments. In their
comparison, the BCC constant is smaller than that of the
$\mathrm{A15}$ quantizer~\cite{SV}.
For quadratic quantization, Ayd{\i}n and Iacobelli establish
asymptotic local optimality of BCC under small periodic
perturbations~\cite{AI}, while Luo and Xie obtain nonlinear
stability estimates uniform over finite-index periods~\cite{LX}.

Among three-dimensional lattices, Barnes and Sloane proved
BCC optimality for $p=2$~\cite{BS}, and Bambah proved it for
the covering radius at fixed covolume~\cite{Bambah}. The latter
is the limiting problem
\[
 \lim_{p\to\infty}E_p(A)^{1/p}
 =|V(A)|^{-1/3}\max_{x\in V(A)}|x| .
\]
Graf and Luschgy asked whether the same lattice minimizes
every finite $p$-th moment
\cite[Conjecture~8.17]{GL}. Bourne, Cagnetti and Roper proved
strict local optimality in lattice shape for every finite
$p\ge1$ and global optimality near $p=2$
\cite[Theorems~1.2--1.3]{BCR}. They also found a local branch
of body-centered tetragonal critical points issuing from the
face-centered cubic lattice at $p=2$
\cite[Theorem~1.5]{BCR}. Theorem~\ref{thm:main} gives the
global minimizer in the Bravais-lattice class for the entire
finite exponent range.

The regular truncated octahedron, the Voronoi cell of BCC,
also minimizes surface area among three-dimensional
parallelohedra of fixed volume~\cite{CN,HS}. The proof by
Hales and Song uses an exact Bernstein certificate for the
final polynomial inequality~\cite[Appendix~A]{HS}.

\paragraph{Outline of the proof.}
\begin{enumerate}
\item \emph{Voronoi moments on a fixed square.}
Section~\ref{sec:geometry} writes the Voronoi cell in Selling
coordinates and subdivides one facet from each opposite pair.
Radial integration over the fifteen tile pyramids and their
opposites gives
\[
 E_p(A)=
 \frac{1}{4(p+3)(\det Q)^{1+p/6}}
 \sum_{T\in\mathcal T}W_T\int_{[-1,1]^2}q_T^{p/2}\,dx\,dy .
\]
Here $Q$ is the Gram matrix of a Selling basis, $q_T$ is the
squared radius on tile $T$, and $W_T$ is its weight. At
$\rho_{ij}=1$ the fifteen quadratics take only two forms.
The BCC equality case is determined by
\[
 \begin{gathered}
  0<\det Q\le16\left(\frac{\sum_{0\le i<j\le3}\rho_{ij}}6\right)^3,\\
  \text{equality}\ \Longleftrightarrow\
  \rho_{01}=\rho_{02}=\rho_{03}=\rho_{12}=\rho_{13}=\rho_{23}.
 \end{gathered}
\]

\item \emph{Polynomial comparison for $1\le p\le35$.}
Section~\ref{sec:bounded} bounds the $p/2$ power of the scaled
ratio of $q_T$ to its BCC value. The lower polynomials are
Taylor polynomials of degrees four and three for $p\in[1,2]$
and $p\in[2,4]$, and weighted tangents of degrees two and three
for $p\in[4,6]$ and $p\in[6,35]$, respectively. Bounds on the
determinant factors then give polynomials in the conorms.
For each of the thirty edge orders up to relabelling, the substitution
\[
 \rho_{\sigma(i)}=X_i+\cdots+X_5,\qquad
 X_0,\ldots,X_5\ge0,\qquad 0\le i\le5
\]
reduces the comparison to coefficientwise nonnegativity.
Numerical Lemmas~\ref{num:taylor-low}--\ref{num:cubic} give
these inequalities by exact interval arithmetic.

\item \emph{Linear comparison for $p\ge35$.}
Section~\ref{sec:large} uses the tangent $s^{p/2}\ge1+\frac p2(s-1)$
and reduces each coefficient to a linear combination of six
weighted reference integrals. Numerical Lemma~\ref{num:tail-finite}
gives nonnegativity on $[35,1028]$. For $p\ge1026$, the signs
follow from estimates at the maximizing corners of the two
reference quadratics. In both comparisons, the determinant
bounds give $E_p(A)>E_p(A_{\mathrm{bcc}})$ when the conorms
are unequal. Equal conorms and $|V(A)|=1$ give the orbit
in~\eqref{eq:minimizers}.

\item \emph{Exact numerical bounds.}
Section~\ref{app:certificate} gives the rational computations for
both comparisons. For bounded exponents, Numerical
Lemmas~\ref{num:taylor-low}--\ref{num:cubic} give coefficient
bounds and counts. Numerical Lemma~\ref{num:tail-forms} computes
the degree-ten coefficient pairs by finite differences. Numerical
Lemma~\ref{num:tail-finite} bounds the corresponding weighted
integrals by Bernstein coefficients.
\end{enumerate}

\section{Selling coordinates and the Voronoi moment}\label{sec:geometry}

We obtain the fifteen fixed-square integrals in~\eqref{eq:moment}
by subdividing the facets of the zonotope~\eqref{eq:zone} and
integrating along rays. The BCC quadratics are given
by~\eqref{eq:q0}. The determinant bound~\eqref{eq:conorm-equality}
identifies the equality case, and~\eqref{eq:sector} gives the
nonnegative variables for the two comparisons.

We use Selling's reduction~\cite{Selling} in the obtuse-superbase
form of Conway and Sloane~\cite[Theorem~8]{CS}.

\begin{lemma}[Selling reduction]\label{lem:selling}
Let $A\in GL(3,\R)$. There exist $t_0,t_1,t_2,t_3\in\Lambda(A)$
such that
\[
 \sum_{i=0}^3t_i=0,\qquad t_i\cdot t_j\le0\quad(i\ne j),
\]
and any three of these vectors form a basis of $\Lambda(A)$.
\end{lemma}

\begin{proof}
Take a basis and the negative of its sum. If $t_i\cdot t_j>0$,
set
\[
 t'_i=-t_i,\qquad t'_j=t_j,\qquad
 t'_k=t_k+t_i\quad(k\notin\{i,j\}),\qquad \sum_{k=0}^3t'_k=0.
\]
For $\{i,j,k,\ell\}=\{0,1,2,3\}$,
\[
 \det(t'_i,t'_k,t'_\ell)=-\det(t_i,t_k,t_\ell).
\]
The relation $\sum_{k=0}^3 t'_k=0$ lets us replace any column of this
basis by $t'_j$, preserving the absolute determinant. Moreover,
\[
 \sum_{k=0}^3\bigl(|t'_k|^2-|t_k|^2\bigr)
 =2|t_i|^2+2t_i\cdot\sum_{k\notin\{i,j\}}t_k
 =-2t_i\cdot t_j<0.
\]
The iteration stays in a bounded subset of $\Lambda(A)$ and
therefore terminates.
\end{proof}

By Lemma~\ref{lem:selling} and~\eqref{eq:basis-invariance}, take
$A=(t_1\mid t_2\mid t_3)$ with $t_0=-t_1-t_2-t_3$ and
$t_i\cdot t_j\le0$ for $i\ne j$. Set $\rho_{ij}=-t_i\cdot t_j$
for $i\ne j$. With $e_0=0$ and
$e_1,e_2,e_3$ the standard coordinate vectors, write $g_{ij}=e_i-e_j$
for the edges $ij$ in $\mathcal E=\{01,02,03,12,13,23\}$. Then
\[
 Q=A^TA=\sum_{e\in\mathcal E}\rho_e g_eg_e^T,\qquad
 \det Q=|V(A)|^2.
\]
In coordinates,
\[
 Q=\begin{pmatrix}
 \rho_{01}+\rho_{12}+\rho_{13} & -\rho_{12} & -\rho_{13}\\
 -\rho_{12} & \rho_{02}+\rho_{12}+\rho_{23} & -\rho_{23}\\
 -\rho_{13} & -\rho_{23} & \rho_{03}+\rho_{13}+\rho_{23}
 \end{pmatrix}.
\]
The inverse relations are
\[
 \rho_{ij}=-Q_{ij}\quad(1\le i<j\le3),\qquad
 \rho_{0i}=\sum_{j=1}^{3}Q_{ij}\quad(1\le i\le3).
\]
Expanding the determinant gives
\begin{equation}\label{eq:tree}
\begin{aligned}
 \det Q={}&\rho_{01}\rho_{02}\rho_{03}\\
 &+\rho_{01}\rho_{02}(\rho_{13}+\rho_{23})
  +\rho_{01}\rho_{03}(\rho_{12}+\rho_{23})
  +\rho_{02}\rho_{03}(\rho_{12}+\rho_{13})\\
 &+(\rho_{01}+\rho_{02}+\rho_{03})
    (\rho_{12}\rho_{13}+\rho_{12}\rho_{23}+\rho_{13}\rho_{23}).
\end{aligned}
\end{equation}

In dimension three, Voronoi cells of lattices are zonotopes
\cite[Sect.~3.2]{SSW}. In the coordinates $z=A^Tx$, the six
generating segments give the Minkowski sum
\begin{equation}\label{eq:zone}
 A^TV(A)=\sum_{e\in\mathcal E}[-\rho_e g_e/2,\rho_e g_e/2],
 \qquad |x|^2=z^TQ^{-1}z .
\end{equation}
For $z$ in the sum and $k\in\Z^3$,
\[
 2z\cdot k\le\sum_{e\in\mathcal E}\rho_e|g_e\cdot k|
 \le\sum_{e\in\mathcal E}\rho_e(g_e\cdot k)^2=k^TQk,
\]
since $g_e\cdot k$ is integral. Thus
\[
 |x-Ak|^2-|x|^2=k^TQk-2z\cdot k\ge0\qquad(k\in\Z^3),
\]
and the sum is contained in $A^TV(A)$.
For three distinct edges, $\det(g_e,g_f,g_h)=0$ on a cycle
and $\det(g_e,g_f,g_h)=\pm1$ on a tree, by expansion at the
leaves. Thus the zonotope volume formula and Cauchy--Binet give
\[
 \begin{aligned}
 \left|\sum_{e\in\mathcal E}[-\rho_e g_e/2,\rho_e g_e/2]\right|
 &=\sum_{e<f<h}\rho_e\rho_f\rho_h\,|\det(g_e,g_f,g_h)|\\
 &=\sum_{e<f<h}\rho_e\rho_f\rho_h\,\det(g_e,g_f,g_h)^2\\
 &=\det Q=|\det A|\,|V(A)|=|A^TV(A)|.
 \end{aligned}
\]
The two full-dimensional compact convex sets are equal by
inclusion and equality of volume.

Suppose first that $\rho_e>0$ for every $e\in\mathcal E$.
For a normal $\eta\ne0$, the supporting face has dimension
\[
 \dim\operatorname{span}\{g_{ij}:\eta\cdot e_i=\eta\cdot e_j\}
 =4-\#\{\eta\cdot e_0,\eta\cdot e_1,\eta\cdot e_2,\eta\cdot e_3\}.
\]
Indeed, within each block of equal values, the differences to one
chosen vertex form a basis for the differences within that block.
The union of these bases is independent by affine independence
of $e_0,e_1,e_2,e_3$. A facet therefore has exactly two blocks.
In each opposite pair take the normal for which the block
containing vertex $3$ has the smaller value, scaled so that
the two values differ by one.

For an edge with $\eta\cdot g_e=0$ the full segment contributes
to the face. For any other edge, $\eta\cdot g_e=\pm1$ and
the contribution is its maximizing endpoint
$\rho_e(\eta\cdot g_e)g_e/2$. The facet center $z_\eta$ and the
supporting plane thus satisfy
\begin{equation}\label{eq:cut-height}
 2z_\eta=\sum_{e\in\mathcal E}\rho_e(\eta\cdot g_e)g_e=Q\eta,
 \qquad
 2\eta\cdot z=\eta^TQ\eta
 \quad\text{on the facet}.
\end{equation}

\begin{lemma}[Tile parametrizations]\label{lem:tiles}
Assume $\rho_e>0$ for every $e\in\mathcal E$. Each
parallelogram facet consists of the single tile
\[
 z_T(x,y)=z_\eta+\tfrac12(x\rho_u g_u+y\rho_v g_v),
\]
where the free-edge pair and normal are
\[
 ((u,v),\eta)\in
 \{((01,23),(0,-1,-1)),\ ((02,13),(-1,0,-1)),\ ((12,03),(1,1,0))\}.
\]
The free edges of each hexagonal facet are $ij,ik,jk$, with
\[
 \begin{aligned}
 ((i,j,k),\eta)\in\{&
 ((1,2,3),(-1,-1,-1)),\ ((0,2,3),(1,0,0)),\\
 &((0,1,3),(0,1,0)),\ ((0,1,2),(0,0,-1))\}.
 \end{aligned}
\]
Each hexagonal facet consists of the three tiles
\begin{align}
 &z_\eta+\tfrac12(x\rho_{ij}g_{ij}+y\rho_{jk}g_{jk}-\rho_{ik}g_{ik}),\notag\\
 &z_\eta+\tfrac12(x\rho_{ij}g_{ij}+y\rho_{ik}g_{ik}+\rho_{jk}g_{jk}),\label{eq:hexagon-tiles}\\
 &z_\eta+\tfrac12(x\rho_{jk}g_{jk}+y\rho_{ik}g_{ik}+\rho_{ij}g_{ij}).\notag
\end{align}
Here $(x,y)\in[-1,1]^2$ and the edges multiplying $x,y$ are
the free edges of the tile.
\end{lemma}

\begin{proof}
The two free segments give the parallelogram. On a hexagonal
facet, write $2(z-z_\eta)=Xg_{ij}+Yg_{jk}$. Since
$g_{ik}=g_{ij}+g_{jk}$, the face consists of the points for
which
\[
 \max\{-\rho_{ik},X-\rho_{ij},Y-\rho_{jk}\}
 \le t\le
 \min\{\rho_{ik},X+\rho_{ij},Y+\rho_{jk}\}
\]
for some $t$, the coefficient of $g_{ik}$ in the segment sum.
Take the left endpoint. The three active bounds give the three
tiles in~\eqref{eq:hexagon-tiles}:
\[
 \begin{aligned}
 t=-\rho_{ik}&\ \Longrightarrow\
      (X,Y)=(x\rho_{ij}-\rho_{ik},y\rho_{jk}-\rho_{ik}),\\
 t=Y-\rho_{jk}&\ \Longrightarrow\
      (X,Y)=(x\rho_{ij}+y\rho_{ik},\rho_{jk}+y\rho_{ik}),\\
 t=X-\rho_{ij}&\ \Longrightarrow\
      (X,Y)=(\rho_{ij}+y\rho_{ik},x\rho_{jk}+y\rho_{ik}).
 \end{aligned}
\]
In each case the other bounds are exactly $|x|,|y|\le1$,
and ties have $|x|=1$ or $|y|=1$.
\end{proof}

Let $\mathcal T$ be the set of $3+4\cdot3=15$ tiles. The
ordered free-edge pair $(u,v)$ identifies each tile uniquely.
On the facet with normal $\eta$, write
\begin{equation}\label{eq:tile-weight}
 q_{u,v}(x,y)=q_T(x,y)=z_T(x,y)^TQ^{-1}z_T(x,y),\qquad
 W_{u,v}=W_T=\rho_u\rho_v\,\eta^TQ\eta .
\end{equation}
Since
$\eta^TQ\eta=\sum_{e\in\mathcal E}\rho_e(\eta\cdot g_e)^2$,
the three parallelogram weights are
\[
\begin{aligned}
 W_{01,23}&=\rho_{01}\rho_{23}(\rho_{02}+\rho_{03}+\rho_{12}+\rho_{13}),\\
 W_{02,13}&=\rho_{02}\rho_{13}(\rho_{01}+\rho_{03}+\rho_{12}+\rho_{23}),\\
 W_{12,03}&=\rho_{12}\rho_{03}(\rho_{01}+\rho_{02}+\rho_{13}+\rho_{23}).
\end{aligned}
\]
The hexagonal weights in the order of~\eqref{eq:hexagon-tiles}
are
\[
\begin{aligned}
 (W_{12,23},W_{12,13},W_{23,13})
   &=(\rho_{01}+\rho_{02}+\rho_{03})
       (\rho_{12}\rho_{23},\rho_{12}\rho_{13},\rho_{23}\rho_{13}),\\
 (W_{02,23},W_{02,03},W_{23,03})
   &=(\rho_{01}+\rho_{12}+\rho_{13})
       (\rho_{02}\rho_{23},\rho_{02}\rho_{03},\rho_{23}\rho_{03}),\\
 (W_{01,13},W_{01,03},W_{13,03})
   &=(\rho_{02}+\rho_{12}+\rho_{23})
       (\rho_{01}\rho_{13},\rho_{01}\rho_{03},\rho_{13}\rho_{03}),\\
 (W_{01,12},W_{01,02},W_{12,02})
   &=(\rho_{03}+\rho_{13}+\rho_{23})
       (\rho_{01}\rho_{12},\rho_{01}\rho_{02},\rho_{12}\rho_{02}).
\end{aligned}
\]
The coordinates in Lemma~\ref{lem:tiles} also give a direct
polynomial evaluation of the squared radius. By~\eqref{eq:cut-height},
\begin{equation}\label{eq:tile-quadratic-centered}
 \begin{split}
  4q_{u,v}\det Q
   ={}&(\eta^TQ\eta)\det Q\\
     &+\bigl(2z_T-2z_\eta\bigr)^T\operatorname{adj}(Q)
             \bigl(2z_T-2z_\eta\bigr).
 \end{split}
\end{equation}
The cross term vanishes because
$z_\eta^TQ^{-1}(z_T-z_\eta)=\eta^T(z_T-z_\eta)/2=0$.
In particular,
\[
 q_T=\tfrac14\eta^TQ\eta+(z_T-z_\eta)^TQ^{-1}(z_T-z_\eta)
     \ge\tfrac14\eta^TQ\eta>0.
\]

\begin{proposition}[Fixed-domain moment identity]\label{prop:moment}
For a reduced basis $A$ with conorms $\rho_e\ge0$, one has
\begin{equation}\label{eq:moment}
 E_p(A)=\frac{1}{4(p+3)|V(A)|^{2+p/3}}
   \sum_{T\in\mathcal T}W_T\int_{[-1,1]^2}q_T(x,y)^{p/2}\,dx\,dy .
\end{equation}
Moreover, $\sum_{T\in\mathcal T}W_T=3\det Q$.
\end{proposition}

\begin{proof}
For the general facet formula for homogeneous integrands,
see Lasserre~\cite{Lasserre}. For positive conorms the pyramid
over $T$ is parametrized by
$z=rz_T(x,y)$, $0\le r\le1$. The free-edge vectors in
Lemma~\ref{lem:tiles} satisfy $g_u\times g_v=\pm\eta$.
By~\eqref{eq:cut-height},
\[
 \big|\det(z_T,\partial_x z_T,\partial_y z_T)\big|
 =\frac{\rho_u\rho_v}{4}|\eta\cdot z_T|=\frac{W_T}{8}.
\]
The radial factor is given by
\[
 dz=\frac{r^2W_T}{8}\,dr\,dx\,dy,\qquad
 z^TQ^{-1}z=r^2q_T(x,y),\qquad \int_0^1r^{p+2}\,dr=\frac1{p+3}.
\]
Including the opposite pyramids gives
\[
 \int_{A^TV(A)}(z^TQ^{-1}z)^{p/2}\,dz
 =\frac1{4(p+3)}\sum_{T\in\mathcal T}W_T
   \int_{[-1,1]^2}q_T(x,y)^{p/2}\,dx\,dy .
\]
Since $|\det A|=|V(A)|$, multiplication by $|V(A)|^{-2-p/3}$
gives~\eqref{eq:moment}. At $p=0$ the unnormalized integral gives
\[
 \det Q=|A^TV(A)|=\frac13\sum_{T\in\mathcal T}W_T.
\]

For vanishing conorms use $\rho_e+\varepsilon$, with basis
$A_\varepsilon=(Q+\varepsilon\sum_{e\in\mathcal E}g_eg_e^T)^{1/2}$.
By~\eqref{eq:zone},
\[
 A^TV(A)\subseteq A_\varepsilon^TV(A_\varepsilon)
 \subseteq A^TV(A)+
 \left\{z:|z|\le\tfrac{\varepsilon}{2}\sum_{e\in\mathcal E}|g_e|\right\}.
\]
Every point outside the compact set $A^TV(A)$ has positive
distance from it. The inclusions therefore give, pointwise,
\[
 \mathbf 1_{A_\varepsilon^TV(A_\varepsilon)}(z)
 \longrightarrow\mathbf 1_{A^TV(A)}(z)\qquad(\varepsilon\downarrow0).
\]
Also $(A_\varepsilon^TA_\varepsilon)^{-1}\to Q^{-1}$ and
$z^T(A_\varepsilon^TA_\varepsilon)^{-1}z\le z^TQ^{-1}z$.
For $0<\varepsilon\le1$ the domains are contained in a fixed
bounded set, on which $(z^TQ^{-1}z)^{p/2}$ is bounded.
Dominated convergence for the integrals and continuity of the
tile expressions give both identities at $\varepsilon=0$.
\end{proof}

\begin{lemma}[BCC reference quadratics]\label{lem:bcc-tiles}
At $\rho=\one$, write $Q_0=Q|_{\rho=\one}$ and
$q_{0,T}=q_T|_{\rho=\one}$. The three square tiles have weight
$4$ and reference quadratic $q_{0,{\rm sq}}$. The twelve rhombus
tiles have weight $3$. On the hexagonal facet with free edges
$ij,ik,jk$, the reference quadratics are
\[
 q_{0,T}(x,y)=
 \begin{cases}
 q_{0,{\rm rh}}(x,y),&(u,v)=(ij,jk),\\
 q_{0,{\rm rh}}(x,-y),&(u,v)\in\{(ij,ik),(jk,ik)\}.
 \end{cases}
\]
Here
\begin{equation}\label{eq:q0}
 \begin{aligned}
 q_{0,{\rm sq}}(x,y)&=1+\frac{x^2+y^2}{8},\\
 q_{0,{\rm rh}}(x,y)&=\frac34+
    \frac{(x-1)^2-(x-1)(y-1)+(y-1)^2}{8}.
 \end{aligned}
\end{equation}
All these quadratics take values in $[3/4,5/4]$ on $[-1,1]^2$.
\end{lemma}

\begin{proof}
At $\rho=\one$ the Gram matrix and its inverse are
\[
 Q_0=\begin{pmatrix}3&-1&-1\\-1&3&-1\\-1&-1&3\end{pmatrix},
 \qquad 4Q_0^{-1}=\begin{pmatrix}2&1&1\\1&2&1\\1&1&2\end{pmatrix},
 \qquad \det Q_0=16.
\]
The columns $(-1,1,1)^T,(1,-1,1)^T,(1,1,-1)^T$ have Gram
matrix $Q_0$. They belong to $2D_3^*$ and generate it, since
both lattices have covolume $4$. By~\eqref{eq:basis-invariance},
the corresponding moment is $E_p(A_{\mathrm{bcc}})$.

At $\rho=\one$ the center contributes
$z_\eta^TQ_0^{-1}z_\eta=\eta^TQ_0\eta/4=W_T/4$ and is
orthogonal to the free directions in the $Q_0^{-1}$ inner
product. The inverse matrix gives
\[
 \begin{aligned}
 g_e^TQ_0^{-1}g_e&=\tfrac12,\qquad
 g_e^TQ_0^{-1}g_f=0 \quad\text{if $e,f$ have disjoint vertices},\\
 g_{ij}^TQ_0^{-1}g_{jk}&=-\tfrac14,\qquad
 g_{ij}^TQ_0^{-1}g_{ik}=g_{jk}^TQ_0^{-1}g_{ik}=\tfrac14
 \quad(i<j<k).
 \end{aligned}
\]
The hexagonal displacements $2(z_T-z_\eta)$ in order are
\[
 (x-1)g_{ij}+(y-1)g_{jk},\qquad
 (x-1)g_{ij}+(y+1)g_{ik},\qquad
 (x-1)g_{jk}+(y+1)g_{ik}.
\]
The listed conorm weights give $W_T=4$ or $3$, so the inner
products yield~\eqref{eq:q0} and its reflections.
The lower bounds follow from~\eqref{eq:q0}, and convexity
gives the upper bounds from the vertex values
\[
 \begin{aligned}
 q_{0,{\rm sq}}(\pm1,\pm1)&=5/4,\\
 q_{0,{\rm rh}}(-1,-1)=q_{0,{\rm rh}}(-1,1)
   =q_{0,{\rm rh}}(1,-1)&=5/4,\qquad
 q_{0,{\rm rh}}(1,1)=3/4.
 \end{aligned}
\]
\end{proof}

Applying~\eqref{eq:moment} at $\rho=\one$ gives
\begin{equation}\label{eq:bcc-value}
 E_p(A_{\mathrm{bcc}})
 =\frac{3}{16(p+3)2^{2p/3}}
 \int_{[-1,1]^2}\bigl(q_{0,{\rm sq}}^{p/2}
                      +3q_{0,{\rm rh}}^{p/2}\bigr)\,dx\,dy .
\end{equation}

For nonnegative integers $a,b$, direct integration gives
\begin{equation}\label{eq:square-monomials}
 \frac14\int_{[-1,1]^2}x^a y^b\,dx\,dy=
 \begin{cases}
 \dfrac{1}{(a+1)(b+1)},&a,b\text{ even},\\
 0,&\text{otherwise}.
 \end{cases}
\end{equation}
For a quadratic $q(x,y)=\sum_{r+s\le2}d_{rs}x^ry^s$, the
multinomial theorem and~\eqref{eq:square-monomials} give
\begin{equation}\label{eq:multinomial-moment}
\begin{aligned}
 &\frac14\int_{[-1,1]^2}x^ay^b q(x,y)^n\,dx\,dy\\
 &\quad=
 \sum_{\substack{\sum_{r+s\le2}m_{rs}=n\\
                  a+\sum_{r+s\le2}r m_{rs}\ \mathrm{even}\\
                  b+\sum_{r+s\le2}s m_{rs}\ \mathrm{even}}}
 \frac{n!\displaystyle\prod_{r+s\le2}
                \bigl(d_{rs}^{m_{rs}}/m_{rs}!\bigr)}
      {\bigl(a+1+\sum_{r+s\le2}r m_{rs}\bigr)
       \bigl(b+1+\sum_{r+s\le2}s m_{rs}\bigr)} .
\end{aligned}
\end{equation}
All indices $r,s,m_{rs}$ and exponents $a,b,n$ are nonnegative
integers. A factor with $m_{rs}=0$ is one.

\begin{numlemma}[First and second moments at BCC]\label{num:moment-one}
One has
\[
 0.46908851641917105<E_1(A_{\mathrm{bcc}})<0.46908851641917107,
 \qquad
 E_2(A_{\mathrm{bcc}})=\frac{19}{64\,2^{1/3}}.
\]
\end{numlemma}

\begin{proof}
Integration of the quadratics in~\eqref{eq:q0} gives
\[
 \int_{[-1,1]^2}q_{0,{\rm sq}}\,dx\,dy=\frac{13}{3},\qquad
 \int_{[-1,1]^2}q_{0,{\rm rh}}\,dx\,dy=\frac{23}{6}.
\]
Thus~\eqref{eq:bcc-value} gives
\[
 E_2(A_{\mathrm{bcc}})
 =\frac{3}{80\cdot2^{4/3}}\left(\frac{13}{3}+3\cdot\frac{23}{6}\right)
 =\frac{19}{64\cdot2^{1/3}}.
\]
Conway and Sloane~\cite[Eq.~(17)]{CSMoments} give the same
moment with the per-coordinate normalization $E_2/3$.

For $q=q_{0,{\rm sq}}$ or $q=q_{0,{\rm rh}}$, the bounds
$|q-1|\le1/4$ and $|\binom{1/2}{j}|\le1$ give
\[
 \left|\frac14\int_{[-1,1]^2}q^{1/2}\,dx\,dy
 -\frac14\sum_{j=0}^{40}\binom{1/2}{j}
      \int_{[-1,1]^2}(q-1)^j\,dx\,dy\right|
 \le\sum_{j=41}^{\infty}4^{-j}=\frac1{3\cdot4^{40}}.
\]
Apply~\eqref{eq:multinomial-moment} to the quadratic $q-1$
with $a=b=0$ and $n=j$. The resulting rational sums and the
error bound give
\[
 \begin{aligned}
 \frac{10\,405\,278\,818\,962\,246\,491}{10^{19}}
 &<\frac14\int_{[-1,1]^2}q_{0,{\rm sq}}^{1/2}\,dx\,dy
 <\frac{10\,405\,278\,818\,962\,246\,492}{10^{19}},\\
 \frac{9\,769\,468\,916\,879\,147\,255}{10^{19}}
 &<\frac14\int_{[-1,1]^2}q_{0,{\rm rh}}^{1/2}\,dx\,dy
 <\frac{9\,769\,468\,916\,879\,147\,256}{10^{19}}.
 \end{aligned}
\]
Cubing rational endpoints gives
$1.587401051968199474<\sqrt[3]{4}<1.587401051968199475$.
Now~\eqref{eq:bcc-value} at $p=1$ reads
\[
 E_1(A_{\mathrm{bcc}})=\frac{3}{64\sqrt[3]{4}}
 \int_{[-1,1]^2}
 \bigl(q_{0,{\rm sq}}^{1/2}+3q_{0,{\rm rh}}^{1/2}\bigr)\,dx\,dy.
\]
The displayed rational bounds give the stated first moment.
\end{proof}

The arithmetic--geometric mean inequality for the eigenvalues
of $Q_0^{-1/2}QQ_0^{-1/2}$ gives
\begin{equation}\label{eq:conorm-equality}
 0<\frac{\det Q}{16}
 \le\left(\frac{\operatorname{tr}(Q_0^{-1/2}QQ_0^{-1/2})}{3}\right)^3
 =\left(\frac{\sum_{e\in\mathcal E}\rho_e}{6}\right)^3.
\end{equation}
The trace is
$\operatorname{tr}(Q_0^{-1}Q)
 =\sum_{e\in\mathcal E}\rho_e g_e^TQ_0^{-1}g_e
 =\tfrac12\sum_{e\in\mathcal E}\rho_e$.
Equality holds exactly when
\[
 Q=\frac{\sum_{e\in\mathcal E}\rho_e}{6}\,Q_0
 \quad\Longleftrightarrow\quad
 \rho_{01}=\rho_{02}=\rho_{03}=\rho_{12}=\rho_{13}=\rho_{23}.
\]
In this case $\Lambda(A)=2(\sum_{e\in\mathcal E}\rho_e/6)^{1/2}
O D_3^*$ for some $O\in O(3)$. At unit covolume,
\[
 1=4\left(\frac{\sum_{e\in\mathcal E}\rho_e}{6}\right)^{3/2},
 \qquad \Lambda(A)=2^{1/3}O D_3^*=OL_{\rm B}
 \quad\text{for some }O\in O(3).
\]
Conversely, every lattice in this orbit has moment
$E_p(A_{\mathrm{bcc}})$ by~\eqref{eq:basis-invariance}.

With $\alpha=p/2$, dividing~\eqref{eq:moment} by~\eqref{eq:bcc-value}
gives
\begin{equation}\label{eq:ratios}
 \frac{E_p(A)}{E_p(A_{\mathrm{bcc}})}
 =\left(\frac{16}{\det Q}\right)^{1+\alpha/3}
 \frac{\displaystyle\sum_{T\in\mathcal T}W_T\int_{[-1,1]^2}
                        q_T^\alpha\,dx\,dy}
      {\displaystyle\int_{[-1,1]^2}
                    \bigl(12q_{0,{\rm sq}}^\alpha+36q_{0,{\rm rh}}^\alpha\bigr)
                    \,dx\,dy}.
\end{equation}

\begin{lemma}[Tile polynomials]\label{lem:tile-polynomial}
For a tile with coordinates $z_T=(z_1,z_2,z_3)^T$,
\begin{equation}\label{eq:N}
\begin{split}
 q_T\det Q={}&
  (\rho_{02}\rho_{03}+\rho_{02}\rho_{23}+\rho_{03}\rho_{23})z_1^2\\
 &+(\rho_{01}\rho_{03}+\rho_{01}\rho_{13}+\rho_{03}\rho_{13})z_2^2\\
 &+(\rho_{01}\rho_{02}+\rho_{01}\rho_{12}+\rho_{02}\rho_{12})z_3^2\\
 &+\rho_{03}\rho_{12}(z_1+z_2)^2
   +\rho_{02}\rho_{13}(z_1+z_3)^2\\
 &+\rho_{01}\rho_{23}(z_2+z_3)^2\\
 &+(\rho_{12}\rho_{13}+\rho_{12}\rho_{23}+\rho_{13}\rho_{23})
                   (z_1+z_2+z_3)^2.
\end{split}
\end{equation}
Equivalently,
\begin{equation}\label{eq:weighted-tile}
 4q_T\det Q=\sum_{\substack{e,f\in\mathcal E\\e<f}}
                   \rho_e\rho_f\det(g_e,g_f,2z_T)^2 .
\end{equation}
For $k\in\mathbb Z_{\ge0}$, the polynomial $W_T(4q_T\det Q)^k$
has integer coefficients, is homogeneous of degree $4k+3$ in
the conorms, and has degree at most $2k$ in $(x,y)$.
\end{lemma}

\begin{proof}
The rank-one determinant identity and Cauchy--Binet give
\[
 \begin{aligned}
 \det(Q+t(2z_T)(2z_T)^T)-\det Q
 &=t(2z_T)^T\operatorname{adj}(Q)(2z_T)\\
 &=t\sum_{e<f}\rho_e\rho_f\det(g_e,g_f,2z_T)^2 .
 \end{aligned}
\]
Since $(2z_T)^T\operatorname{adj}(Q)(2z_T)=4q_T\det Q$,
this gives~\eqref{eq:weighted-tile}. For each pair of edges,
\[
 \det(g_e,g_f,2z_T)^2=4\bigl((g_e\times g_f)\cdot z_T\bigr)^2.
\]
The seven possible cross products, up to sign, are
\[
 e_1,\ e_2,\ e_3,\ e_1+e_2,\ e_1+e_3,\ e_2+e_3,\ e_1+e_2+e_3.
\]
Grouping the fifteen pairs by their cross products gives~\eqref{eq:N}.
The tile parametrizations give $2z_T$ with integer coefficients,
linear homogeneous in the conorms and affine in $(x,y)$.
The weight in~\eqref{eq:tile-weight} is cubic with integer
coefficients, giving the claimed degrees.
\end{proof}

Order the edges as $01,02,03,12,13,23$. If
$\sigma:\{0,\ldots,5\}\to\mathcal E$ is a
bijection with $\rho_{\sigma(0)}\ge\cdots\ge\rho_{\sigma(5)}$, write
\begin{equation}\label{eq:sector}
 \rho_{\sigma(i)}=\sum_{j=i}^{5}X_j,\qquad
 X_i\ge0,\qquad 0\le i\le5.
\end{equation}
Thus $X_i=\rho_{\sigma(i)}-\rho_{\sigma(i+1)}$ for $0\le i<5$
and $X_5=\rho_{\sigma(5)}$. There are $6!/4!=30$ orderings
modulo vertex relabellings, since a vertex permutation fixing
each edge is the identity.
In each orbit take the lexicographically least order. Relabelling
the superbase preserves $E_p(A)$ by~\eqref{eq:basis-invariance}.
The closures of the thirty sectors include the cases with ties.

For $\mu\in\mathbb Z_{\ge0}^6$, expanding the powers of
$X_0,\ldots,X_5$ in~\eqref{eq:sector} gives
\begin{equation}\label{eq:sector-coeff}
 \prod_{i=0}^5\left(\sum_{j=i}^5X_j\right)^{\mu_i}
 =
 \sum_{\substack{\nu\in\mathbb Z_{\ge0}^6\\
                  \sum_{i=0}^5\nu_i=\sum_{i=0}^5\mu_i}}
 \left(
 \prod_{j=0}^5
 \binom{\displaystyle\sum_{i=0}^j\mu_i-\sum_{i=0}^{j-1}\nu_i}{\nu_j}
 \right)\prod_{i=0}^5 X_i^{\nu_i}.
\end{equation}
The second sum in the upper entry is empty at $j=0$. A factor
vanishes if its upper entry is smaller than its lower entry.
Indeed, at step $j$ there are
$\sum_{i=0}^j\mu_i-\sum_{i=0}^{j-1}\nu_i$ powers of
$\sum_{i=j}^5 X_i$ from which to select $X_j^{\nu_j}$.

For homogeneous $P$ of degree $n\ge1$, the terms of degree
at most one in $X_0,\ldots,X_4$ in $P(\rho(X))$ are
\[
 P(\one)X_5^n+
 X_5^{n-1}\sum_{i=0}^4 X_i
              \sum_{j=0}^{i}\partial_{\rho_{\sigma(j)}}P(\one).
\]
Indeed, with $X_0=\cdots=X_4=0$ and $X_5=1$,
\eqref{eq:sector} gives $\rho=\one$ and the derivative in $X_i$
is the sum of the derivatives in $\rho_{\sigma(0)},\ldots,
\rho_{\sigma(i)}$. We call these six monomials the first-order
positions.

\Needspace{11\baselineskip}

\section{Comparison for bounded exponents}\label{sec:bounded}

We prove BCC optimality and uniqueness for $1\le p\le35$.
A lower bound for the power of the scaled squared radius is combined
with upper bounds for the fractional determinant powers. The resulting
coefficient inequalities are quartic Taylor for $1\le p\le2$,
cubic Taylor for $2\le p\le4$, weighted quadratic for $4\le p\le6$,
and weighted cubic for $6\le p\le35$.

\begin{proposition}\label{prop:bounded}
Theorem~\ref{thm:main} holds for $1\le p\le35$.
\end{proposition}

\begin{lemma}
Let $0<\beta<1$ and $0<t\le1$. Then
\begin{align}
 t^\beta&\le\frac{(1-\beta)+(1+\beta)t}{(1+\beta)+(1-\beta)t},
 \label{eq:pade}\\
 t^{j/3}&\le 1-\frac j3(1-t),\qquad j=1,2,3.
 \label{eq:tangents}
\end{align}
For $t<1$ the inequality in~\eqref{eq:pade} is strict, as is
\eqref{eq:tangents} for $j=1,2$.
\end{lemma}

\begin{proof}
Write $t=e^{-2u}$. Strict concavity of $\tanh$ on $[0,\infty)$ gives
\[
 \frac{1-t^\beta}{1+t^\beta}
 =\tanh(\beta u)\ge\beta\tanh u=\beta\frac{1-t}{1+t}.
\]
The inequality is strict for $u>0$ and rearranges to~\eqref{eq:pade}.
For~\eqref{eq:tangents}, use concavity of $t^{j/3}$, which is
strict for $j=1,2$.
\end{proof}

\begin{lemma}
For $r\ge0$ and $1/2\le\alpha\le2$, one has
\begin{equation}\label{eq:taylor}
 r^\alpha\ge
 \sum_{j=0}^3\binom{\alpha}{j}(r-1)^j+
 \begin{cases}
 \displaystyle\binom{\alpha-1}{3}(r-1)^4,&1/2\le\alpha\le1,\\[2pt]
 0,&1\le\alpha\le2.
 \end{cases}
\end{equation}
\end{lemma}

\begin{proof}
For $r>0$, $r\ne1$, the integral remainder gives
\[
 \frac{r^\alpha-\sum_{j=0}^3\binom{\alpha}{j}(r-1)^j}{(r-1)^4}
 =\frac{\alpha(\alpha-1)(\alpha-2)(\alpha-3)}{3!}
 \int_0^1(1-\tau)^3(1+\tau(r-1))^{\alpha-4}\,d\tau .
\]
For $1\le\alpha\le2$ the right-hand side is nonnegative.
For $1/2\le\alpha<1$ its prefactor is negative, and
\[
 (1+\tau(r-1))^{\alpha-4}\le(1-\tau)^{\alpha-4},
 \qquad 0\le\tau<1 .
\]
Thus the quotient is bounded below by
\[
 \frac{\alpha(\alpha-1)(\alpha-2)(\alpha-3)}{3!}
 \int_0^1(1-\tau)^{\alpha-1}\,d\tau
 =\binom{\alpha-1}{3}.
\]
Continuity gives~\eqref{eq:taylor} also at $r=0,1$.
\end{proof}

For $\alpha=p/2$ and $k\in\mathbb Z_{\ge0}$, let
\[
 I_k=\frac14\sum_{T\in\mathcal T}W_T\int_{[-1,1]^2}
        (4q_T\det Q)^k q_{0,T}^{\alpha-k}\,dx\,dy .
\]
By~\eqref{eq:q0}, the parallelogram terms use $q_{0,{\rm sq}}$.
The hexagon terms use $q_{0,{\rm rh}}$ after $y\mapsto-y$
on the tiles with pairs $(ij,i\ell)$ and $(j\ell,i\ell)$. Thus
\begin{equation}\label{eq:ik-explicit}
\begin{aligned}
 \frac{4I_k}{(4\det Q)^k}
 ={}&\int_{[-1,1]^2}q_{0,{\rm sq}}^{\alpha-k}
       \bigl[W_{01,23}q_{01,23}(x,y)^k
             +W_{02,13}q_{02,13}(x,y)^k\\
 &\hspace{34mm}{}+W_{12,03}q_{12,03}(x,y)^k\bigr]\,dx\,dy\\
 &+\sum_{\substack{(i,j,\ell)\in\{(1,2,3),(0,2,3),\\
                              (0,1,3),(0,1,2)\}}}
       \int_{[-1,1]^2}q_{0,{\rm rh}}^{\alpha-k}
       \bigl[W_{ij,j\ell}q_{ij,j\ell}(x,y)^k\\
 &\hspace{34mm}{}+W_{ij,i\ell}q_{ij,i\ell}(x,-y)^k
                +W_{j\ell,i\ell}q_{j\ell,i\ell}(x,-y)^k\bigr]\,dx\,dy .
\end{aligned}
\end{equation}
Taking $k=0$ and using $\sum_{T\in\mathcal T}W_T=3\det Q$ gives
\begin{equation}\label{eq:i0-explicit}
\begin{aligned}
 4I_0=\int_{[-1,1]^2}\Bigl[
     & (W_{01,23}+W_{02,13}+W_{12,03})
                 (q_{0,{\rm sq}}^\alpha-q_{0,{\rm rh}}^\alpha)\\
     &+3(\det Q)q_{0,{\rm rh}}^\alpha\Bigr]\,dx\,dy .
\end{aligned}
\end{equation}
With $4q_T\det Q$ given by~\eqref{eq:N}, expand
\[
 W_T(4q_T\det Q)^k=
 \sum_{\substack{\nu\in\mathbb Z_{\ge0}^6,\
                 \sum_{e\in\mathcal E}\nu_e=4k+3\\
                 a,b\in\mathbb Z_{\ge0},\ a+b\le2k}}
 d_{\nu,a,b}(T,k)\left(\prod_{e\in\mathcal E}\rho_e^{\nu_e}\right)x^ay^b.
\]
Termwise integration yields the conorm polynomial
\begin{equation}\label{eq:ik-coeff}
\begin{split}
 I_k={}&\frac14\sum_{T\in\mathcal T}
 \sum_{\substack{\nu\in\mathbb Z_{\ge0}^6,\ \sum_{e\in\mathcal E}\nu_e=4k+3\\
                 a,b\in\mathbb Z_{\ge0},\ a+b\le2k}}
 d_{\nu,a,b}(T,k)\Bigl(\prod_{e\in\mathcal E}\rho_e^{\nu_e}\Bigr)
       \int_{[-1,1]^2}x^ay^bq_{0,T}^{\alpha-k}\,dx\,dy .
\end{split}
\end{equation}

Write the right-hand side of~\eqref{eq:taylor} as
$\sum_{k=0}^K c_k(\alpha)r^k$, with $K=4$ on $[1/2,1]$ and $K=3$
on $[1,2]$. The two polynomials agree at $\alpha=1$. Expanding
the powers of $r-1$ and using
$\sum_{j=0}^m(-1)^j\binom{x}{j}=(-1)^m\binom{x-1}{m}$ gives
\[
 c_k(\alpha)=
 \begin{cases}
 (-1)^{3-k}\binom{\alpha}{k}\binom{\alpha-k-1}{3-k},
   &K=3,\quad 0\le k\le3,\\
 (-1)^k\left[\binom4k\binom{\alpha-1}{3}
              -\binom{\alpha}{k}\binom{\alpha-k-1}{3-k}\right],
   &K=4,\quad 0\le k\le4,
 \end{cases}
\]
where $\binom{x}{-1}=0$. In both cases $\sum_{k=0}^K c_k(\alpha)=1$.

The Taylor comparison uses the polynomial
\begin{equation}\label{eq:lowpoly}
\begin{split}
 P_\alpha^{\rm T}={}&
 \left[\left(1+\frac{\alpha}{3}\right)
       \Bigl(\sum_{e\in\mathcal E}\rho_e\Bigr)^3+\tfrac92(3-\alpha)\det Q\right]\\
 &\qquad{}\times\sum_{k=0}^K3^kc_k(\alpha)
       \Bigl(2\det Q\sum_{e\in\mathcal E}\rho_e\Bigr)^{K-k}I_k\\
 &-\frac{\det Q}{16}\Bigl(2\det Q\sum_{e\in\mathcal E}\rho_e\Bigr)^K
   \left[\left(1-\frac{\alpha}{3}\right)
       \Bigl(\sum_{e\in\mathcal E}\rho_e\Bigr)^3+\tfrac92(3+\alpha)\det Q\right]\\
 &\qquad{}\times
   \int_{[-1,1]^2}\bigl(3q_{0,{\rm sq}}^\alpha
                          +9q_{0,{\rm rh}}^\alpha\bigr)\,dx\,dy .
\end{split}
\end{equation}
The quadratic and cubic comparisons use
\begin{align}
 P_\alpha^{(2)}
   ={}&864(\alpha-1)\Bigl(\sum_{e\in\mathcal E}\rho_e\Bigr)^3 I_2
       \notag\\
     &-96(\alpha-2)\det Q
       \Bigl(\sum_{e\in\mathcal E}\rho_e\Bigr)
       \left[4\Bigl(\sum_{e\in\mathcal E}\rho_e\Bigr)^3+27\det Q\right]I_1
       \notag\\
     &-\det Q\Bigl(2\det Q\sum_{e\in\mathcal E}\rho_e\Bigr)^2
       \left[2\Bigl(\sum_{e\in\mathcal E}\rho_e\Bigr)^3+54\det Q\right]
       \notag\\
     &\qquad{}\times\int_{[-1,1]^2}
          \bigl(3q_{0,{\rm sq}}^\alpha+9q_{0,{\rm rh}}^\alpha\bigr)\,dx\,dy ,
       \label{eq:quadpoly}\\
 P_\alpha^{(3)}
   ={}&288(\alpha-2)\Bigl(\sum_{e\in\mathcal E}\rho_e\Bigr)^3 I_3
       \notag\\
     &-32(\alpha-3)\det Q
       \Bigl(\sum_{e\in\mathcal E}\rho_e\Bigr)
       \left[4\Bigl(\sum_{e\in\mathcal E}\rho_e\Bigr)^3+27\det Q\right]I_2
       \notag\\
     &-9(\det Q)^2\Bigl(2\det Q\sum_{e\in\mathcal E}\rho_e\Bigr)^3
       \int_{[-1,1]^2}\bigl(3q_{0,{\rm sq}}^\alpha
                              +9q_{0,{\rm rh}}^\alpha\bigr)\,dx\,dy .
       \label{eq:cubicpoly}
\end{align}

The moments $I_k$ are invariant under vertex relabelling
by~\eqref{eq:ik-relabel}, as are $\det Q$ and
$\sum_{e\in\mathcal E}\rho_e$. Each comparison polynomial
is therefore invariant. Every conorm tuple is relabelled into
one of the thirty sectors in~\eqref{eq:sector}, where
$X_i\ge0$ for all $i$. Nonnegative coefficients in these sectors
imply nonnegativity for every conorm tuple.

\begin{proof}[Proof of Proposition~\ref{prop:bounded}]
By~\eqref{eq:conorm-equality},
\[
 0<t=\frac{27\det Q}{2\bigl(\sum_{e\in\mathcal E}\rho_e\bigr)^3}\le1.
\]
For each tile,
\begin{equation}\label{eq:taylor-argument}
 \frac{6q_T}{q_{0,T}\sum_{e\in\mathcal E}\rho_e}
 =t^{1/3}\left(\frac{16}{\det Q}\right)^{1/3}\frac{q_T}{q_{0,T}}
 =\frac{3(4q_T\det Q)}{2q_{0,T}(\det Q)\sum_{e\in\mathcal E}\rho_e} .
\end{equation}
In particular,
\begin{equation}\label{eq:scaled-moments}
 \frac14\sum_{T\in\mathcal T}W_T
 \int_{[-1,1]^2} q_{0,T}^\alpha
       \left(\frac{6q_T}{q_{0,T}\sum_{e\in\mathcal E}\rho_e}\right)^k\,dx\,dy
 =\left(\frac{3}{2\det Q\sum_{e\in\mathcal E}\rho_e}\right)^k I_k .
\end{equation}

\emph{Case 1. $1\le p\le4$.}
By~\eqref{eq:taylor}, \eqref{eq:scaled-moments} and~\eqref{eq:ratios},
\[
 \begin{aligned}
 &\frac{\det Q}{16}\,t^{\alpha/3}\frac{E_p(A)}{E_p(A_{\mathrm{bcc}})}
   \int_{[-1,1]^2}\bigl(3q_{0,{\rm sq}}^\alpha+9q_{0,{\rm rh}}^\alpha\bigr)\,dx\,dy\\
 &\qquad\ge\sum_{k=0}^K c_k(\alpha)
       \left(\frac{3}{2\det Q\sum_{e\in\mathcal E}\rho_e}\right)^k I_k\\
 &\qquad\ge\frac{\det Q}{16}
       \frac{(1-\alpha/3)+(1+\alpha/3)t}
            {(1+\alpha/3)+(1-\alpha/3)t}
       \int_{[-1,1]^2}\bigl(3q_{0,{\rm sq}}^\alpha+9q_{0,{\rm rh}}^\alpha\bigr)\,dx\,dy.
 \end{aligned}
\]
The last inequality is $P_\alpha^{\rm T}\ge0$ from Numerical
Lemmas~\ref{num:taylor-low} and~\ref{num:taylor-high}, divided by
\[
 \Bigl(2\det Q\sum_{e\in\mathcal E}\rho_e\Bigr)^K
 \Bigl(\sum_{e\in\mathcal E}\rho_e\Bigr)^3
 \bigl[(1+\alpha/3)+(1-\alpha/3)t\bigr]>0.
\]
By~\eqref{eq:pade},
\[
 \frac{E_p(A)}{E_p(A_{\mathrm{bcc}})}
 \ge\frac{(1-\alpha/3)+(1+\alpha/3)t}
          {t^{\alpha/3}\bigl[(1+\alpha/3)+(1-\alpha/3)t\bigr]}
 \ge1,
\]
with strict inequality for $t<1$.

\emph{Case 2. $4\le p\le35$.}
Take $k=2$ for $4\le p\le6$ and $k=3$ for $6<p\le35$.
The tangent to $s^{\alpha-k+1}$ at one gives
\begin{equation}\label{eq:power}
 s^\alpha\ge(\alpha-k+1)s^k-(\alpha-k)s^{k-1},\qquad s\ge0.
\end{equation}
Use~\eqref{eq:power} with $s=(16/\det Q)^{1/3}q_T/q_{0,T}$
in the tile integrals. Equations~\eqref{eq:ratios},
\eqref{eq:taylor-argument} and~\eqref{eq:scaled-moments} give
\[
\begin{aligned}
 &\frac{\det Q}{16}\Bigl(2\det Q\sum_{e\in\mathcal E}\rho_e\Bigr)^k
       t^{k/3}\left(\frac{E_p(A)}{E_p(A_{\mathrm{bcc}})}-1\right)
       \int_{[-1,1]^2}\bigl(3q_{0,{\rm sq}}^\alpha+9q_{0,{\rm rh}}^\alpha\bigr)\,dx\,dy\\
 &\quad\ge 3^k(\alpha-k+1)I_k
      -3^{k-1}(\alpha-k)\Bigl(2\det Q\sum_{e\in\mathcal E}\rho_e\Bigr)t^{1/3}I_{k-1}\\
 &\qquad\quad-\frac{\det Q}{16}\Bigl(2\det Q\sum_{e\in\mathcal E}\rho_e\Bigr)^k t^{k/3}
       \int_{[-1,1]^2}\bigl(3q_{0,{\rm sq}}^\alpha+9q_{0,{\rm rh}}^\alpha\bigr)\,dx\,dy\\
 &\quad\ge
 \frac{1}{96\bigl(\sum_{e\in\mathcal E}\rho_e\bigr)^3}
 \begin{cases}
 P_\alpha^{(2)},&k=2,\\
 9P_\alpha^{(3)},&k=3
 \end{cases}
 \ \ge0.
\end{aligned}
\]
The second inequality uses~\eqref{eq:tangents} and the identities
\[
 \begin{aligned}
  6\Bigl(\sum_{e\in\mathcal E}\rho_e\Bigr)^3(2/3+t/3)&=4\Bigl(\sum_{e\in\mathcal E}\rho_e\Bigr)^3+27\det Q,\\
  6\Bigl(\sum_{e\in\mathcal E}\rho_e\Bigr)^3(1/3+2t/3)
       &=2\Bigl(\sum_{e\in\mathcal E}\rho_e\Bigr)^3+54\det Q,\\
  6\Bigl(\sum_{e\in\mathcal E}\rho_e\Bigr)^3t&=81\det Q.
 \end{aligned}
\]
The last inequality is Numerical Lemma~\ref{num:quadratic}
for $k=2$ and Numerical Lemma~\ref{num:cubic} for $k=3$.

For $t<1$, the error from the tangent bounds in
\eqref{eq:tangents} includes a positive multiple of
\[
 \begin{cases}
 \displaystyle
 \left(\tfrac13+\tfrac23t-t^{2/3}\right)
 \int_{[-1,1]^2}
     (3q_{0,{\rm sq}}^\alpha+9q_{0,{\rm rh}}^\alpha)\,dx\,dy,
       & k=2,\\[3pt]
 \displaystyle
 (\alpha-3)\left(\tfrac23+\tfrac13t-t^{1/3}\right)I_2,
       & k=3.
 \end{cases}
\]
The integral is positive. Also $\alpha>3$ when $k=3$, and
$I_2>0$ since $\sum_{T\in\mathcal T}W_T=3\det Q>0$
and $q_T\ge\eta^TQ\eta/4>0$ by~\eqref{eq:tile-quadratic-centered}.
Thus the comparison is strict for $t<1$. By
\eqref{eq:conorm-equality}, $t=1$ gives precisely the
unit-covolume orbit in~\eqref{eq:minimizers}.
\end{proof}

At $\rho=\one$, the values in the comparison polynomials are
\[
 \det Q=16,\qquad \sum_{e\in\mathcal E}\rho_e=6,\qquad
 I_k(\one)=64^k\int_{[-1,1]^2}
       \bigl(3q_{0,{\rm sq}}^\alpha+9q_{0,{\rm rh}}^\alpha\bigr)\,dx\,dy .
\]
In~\eqref{eq:lowpoly}, both bracketed factors equal $432$,
and $2\det Q\sum_{e\in\mathcal E}\rho_e=192=3\cdot64$.
Consequently,
\[
 \begin{aligned}
 P_\alpha^{\rm T}(\one)
 &=432\cdot192^K I_0(\one)\left(\sum_{k=0}^K c_k(\alpha)-1\right)=0,\\
 P_\alpha^{(2)}(\one)
 &=864\cdot6^3I_2(\one)\bigl((\alpha-1)-(\alpha-2)-1\bigr)=0,\\
 P_\alpha^{(3)}(\one)
 &=288\cdot6^3I_3(\one)\bigl((\alpha-2)-(\alpha-3)-1\bigr)=0.
 \end{aligned}
\]
Each polynomial is homogeneous and, by~\eqref{eq:ik-relabel},
invariant under relabelling the superbase.
Since the vertex permutations act transitively on the six edges,
Euler's identity gives
\[
 \partial_{\rho_e}P(\one)
 =\frac16\sum_{f\in\mathcal E}\partial_{\rho_f}P(\one)
 =\frac{\deg P}{6}P(\one)=0 .
\]
Hence the six first-order coefficients in each sector vanish.

\section{Comparison for large exponents}\label{sec:large}

For $p\ge35$, the tangent inequality for the normalized squared
radius reduces the comparison to coefficient inequalities for a
degree-ten polynomial. The six BCC integrals in its coefficients
are controlled by Bernstein
bounds for $35\le p\le1028$ and by an estimate at the maximizing
corners for $p\ge1026$.

\begin{proposition}\label{prop:large}
For $35\le p<\infty$ and $A\in GL(3,\R)$,
\[
 E_p(A)\ge E_p(A_{\mathrm{bcc}}).
\]
Among lattices of unit covolume, equality holds exactly on the
orbit in~\eqref{eq:minimizers}.
\end{proposition}

Let $\alpha=p/2$. By~\eqref{eq:q0}, $3/5\le 4q_{0,T}/5\le1$
on $[-1,1]^2$.

\begin{lemma}[Polynomial comparison]\label{lem:tail-reduction}
Let $\alpha\ge35/2$. For a reduced basis matrix $A$ with conorms
$\rho$, define
\begin{equation}\label{eq:tailpoly}
\begin{split}
 P_\alpha^{(1)}
  ={}&384\alpha\left[4\Bigl(\sum_{e\in\mathcal E}\rho_e\Bigr)^3+27\det Q\right]\\
    &\qquad{}\times\sum_{T\in\mathcal T} W_T\int_{[-1,1]^2}
       q_T(\det Q)\left(\frac45q_{0,T}\right)^{\alpha-1}\,dx\,dy\\
    &-5\det Q\Bigl(\sum_{e\in\mathcal E}\rho_e\Bigr)
          \left[2\Bigl(\sum_{e\in\mathcal E}\rho_e\Bigr)^3+54\det Q\right]\\
    &\qquad{}\times\biggl[
       16(\alpha-1)\sum_{T\in\mathcal T}W_T\int_{[-1,1]^2}
                    \left(\frac45q_{0,T}\right)^\alpha\,dx\,dy\\
    &\hspace{28mm}
       +(\det Q)(4/5)^\alpha\int_{[-1,1]^2}
               (12q_{0,{\rm sq}}^\alpha+36q_{0,{\rm rh}}^\alpha)\,dx\,dy\biggr].
\end{split}
\end{equation}
If $P_\alpha^{(1)}\ge0$, then $E_p(A)\ge E_p(A_{\mathrm{bcc}})$.
The inequality is strict when the conorms are unequal.
\end{lemma}

\begin{proof}
Apply $s^\alpha\ge\alpha s+1-\alpha$ with
$s=(16/\det Q)^{1/3}q_T/q_{0,T}$ and multiply by $q_{0,T}^\alpha$.
This gives
\[
 \left(\frac{16}{\det Q}\right)^{\alpha/3}q_T^\alpha
 \ge\alpha\left(\frac{16}{\det Q}\right)^{1/3}q_Tq_{0,T}^{\alpha-1}
       -(\alpha-1)q_{0,T}^\alpha .
\]
For $t=27\det Q/(2(\sum_{e\in\mathcal E}\rho_e)^3)\in(0,1]$,
\eqref{eq:pade} with $\beta=1/3$ gives
\[
 \left(\frac{16}{\det Q}\right)^{1/3}
 =\frac{6t^{-1/3}}{\sum_{e\in\mathcal E}\rho_e}
 \ge \frac{6}{\sum_{e\in\mathcal E}\rho_e}
      \frac{4(\sum_{e\in\mathcal E}\rho_e)^3+27\det Q}
           {2(\sum_{e\in\mathcal E}\rho_e)^3+54\det Q} .
\]
Integrate the tangent bound against $W_T$ and sum over the tiles.
Using~\eqref{eq:ratios} and then substituting the Pad\'e bound gives
\[
\begin{aligned}
 &\left(\frac{E_p(A)}{E_p(A_{\mathrm{bcc}})}-1\right)
       (\det Q)\int_{[-1,1]^2}
           \bigl(3q_{0,{\rm sq}}^\alpha+9q_{0,{\rm rh}}^\alpha\bigr)\,dx\,dy\\
 &\qquad\ge
  \frac{4\alpha}{\det Q}\left(\frac{16}{\det Q}\right)^{1/3} I_1
  -16(\alpha-1)I_0
  -(\det Q)\int_{[-1,1]^2}
           \bigl(3q_{0,{\rm sq}}^\alpha+9q_{0,{\rm rh}}^\alpha\bigr)\,dx\,dy\\
 &\qquad\ge
 \frac{P_\alpha^{(1)}}
 {20(4/5)^\alpha(\det Q)(\sum_{e\in\mathcal E}\rho_e)
                     \bigl[2(\sum_{e\in\mathcal E}\rho_e)^3+54\det Q\bigr]} .
\end{aligned}
\]
The numerator in the last line is exactly~\eqref{eq:tailpoly}.
For unequal conorms the bound from~\eqref{eq:pade} is strict, and
$I_1>0$ since $\sum_{T\in\mathcal T}W_T=3\det Q>0$ and
$q_T\ge\eta^TQ\eta/4>0$ by~\eqref{eq:tile-quadratic-centered}.
\end{proof}

Write the six reference integrals as
\[
 F_\alpha=\frac{(4/5)^{\alpha-1}}4
 \left(\int_{[-1,1]^2}(1,x^2)q_{0,{\rm sq}}^{\alpha-1}\,dx\,dy,\
       \int_{[-1,1]^2}(1,x,x^2,xy)q_{0,{\rm rh}}^{\alpha-1}\,dx\,dy
 \right).
\]
For $R(x,y)=\sum_{i+j\le2}r_{ij}x^iy^j$, define the vector
$\Phi_T(R)$ by
\begin{equation}\label{eq:tail-weight-vectors}
 \Phi_T(R)=
 \begin{cases}
 (r_{00},r_{20}+r_{02},0,0,0,0),&q_{0,T}=q_{0,{\rm sq}},\\
 (0,0,r_{00},r_{10}+r_{01},r_{20}+r_{02},r_{11}),
 &q_{0,T}=q_{0,{\rm rh}}(x,y),\\
 (0,0,r_{00},r_{10}-r_{01},r_{20}+r_{02},-r_{11}),
 &q_{0,T}=q_{0,{\rm rh}}(x,-y).
 \end{cases}
\end{equation}
The symmetry in $x,y$, evenness on the square, and the change
$y\mapsto-y$ on the reflected rhombus give
\[
 \Phi_T(R)\cdot F_\alpha
 =\frac14\int_{[-1,1]^2}R(4q_{0,T}/5)^{\alpha-1}\,dx\,dy .
\]

A vertex permutation $\pi$ acts on the dual coordinates by
$U_\pi e_i=e_{\pi(i)}-e_{\pi(0)}$ for $i=1,2,3$. With primes denoting
the relabelled data,
\[
 U_\pi g_{ij}=e_{\pi(i)}-e_{\pi(j)},\qquad
 Q'=U_\pi Q U_\pi^T,\qquad
 \det Q'=\det Q,\qquad \sum_{e\in\mathcal E}\rho'_e=\sum_{e\in\mathcal E}\rho_e .
\]
Here $U_\pi^{-1}$ is the map associated to $\pi^{-1}$, so $\det U_\pi=\pm1$.
The relabelled facet representative has normal
$\eta'=\pm U_\pi^{-T}\eta$, hence $\eta'^TQ'\eta'=\eta^TQ\eta$.
Match the free-edge pair, and on each hexagon also the remaining
edge that determines the tile's shift in~\eqref{eq:hexagon-tiles}.
For each matched pair of tiles there is a signed permutation
$(x,y)\mapsto(x',y')$ for which
\[
 z'_{T'}(x',y')-z'_{\eta'}
  =\pm U_\pi(z_T(x,y)-z_\eta).
\]
On a hexagon choose the overall sign to match the fixed-edge
shift, and the two variable signs to match the free directions.
By~\eqref{eq:tile-weight} and~\eqref{eq:tile-quadratic-centered},
\[
 W'_{T'}=W_T,\qquad
 q'_{T'}(x',y')=q_T(x,y),\qquad
 q'_{0,T'}(x',y')=q_{0,T}(x,y).
\]
The last identity is the preceding one at $\rho=\one$. The
changes of variables preserve $[-1,1]^2$ and its area, so
\begin{equation}\label{eq:ik-relabel}
 I_k(\rho')=I_k(\rho)\qquad(k\ge0).
\end{equation}
Consequently the sum of $P_\alpha^{(1)}$ over the twenty-four
relabellings is $P_\alpha^{\rm sym}=24P_\alpha^{(1)}$.

For the weights $16q_{0,T}$ in~\eqref{eq:tail-weight-vectors},
the three square tiles give $(16,4,0,0,0,0)$ and the twelve
rhombus tiles give $(0,0,14,-4,4,-2)$. Thus
\[
 \begin{aligned}
 U(\rho)
 &=\sum_{T\in\mathcal T}W_T\Phi_T(16q_{0,T})\\
 &=(W_{01,23}+W_{02,13}+W_{12,03})(16,4,-14,4,-4,2)\\
 &\quad{}+3(\det Q)(0,0,14,-4,4,-2),\\
 U(\one)&=(192,48,504,-144,144,-72).
 \end{aligned}
\]
The first tile sum in~\eqref{eq:tailpoly} is
$\sum_{T\in\mathcal T}W_T\bigl(\Phi_T(4q_T\det Q)\cdot F_\alpha\bigr)$.
Also
\[
 U(\rho)\cdot F_\alpha
   =5\sum_{T\in\mathcal T}W_T
           \int_{[-1,1]^2}(4q_{0,T}/5)^\alpha\,dx\,dy .
\]
The last bracket in~\eqref{eq:tailpoly} is therefore
$\bigl[16(\alpha-1)U(\rho)+(\det Q)U(\one)\bigr]\cdot F_\alpha/5$,
so
\begin{equation}\label{eq:tail-sum}
\begin{aligned}
 P_\alpha^{\rm sym}(\rho)
 ={}&9216\alpha\left[4\Bigl(\sum_{e\in\mathcal E}\rho_e\Bigr)^3
                           +27\det Q\right]\\
 &\qquad{}\times\sum_{T\in\mathcal T}W_T
                    \bigl(\Phi_T(4q_T\det Q)\cdot F_\alpha\bigr)\\
 &-24(\det Q)\Bigl(\sum_{e\in\mathcal E}\rho_e\Bigr)
        \left[2\Bigl(\sum_{e\in\mathcal E}\rho_e\Bigr)^3+54\det Q\right]\\
 &\qquad{}\times
       \bigl(\left[16(\alpha-1)U(\rho)+(\det Q)U(\one)\right]
             \cdot F_\alpha\bigr).
\end{aligned}
\end{equation}
Numerical Lemma~\ref{num:tail-forms} gives the vectors
$C_{\sigma,\nu}$ and $B_{\sigma,\nu}$ by the direct sum
\eqref{eq:tail-sector-pair}. The coefficient of
$X_0^{\nu_0}\cdots X_5^{\nu_5}$ in~\eqref{eq:tail-sum} after
$\rho_{\sigma(i)}=X_i+\cdots+X_5$ is
$(C_{\sigma,\nu}+\alpha B_{\sigma,\nu})\cdot F_\alpha$.
The set $\mathcal C$ consists of the distinct pairs obtained
by dividing each nonzero pair in~\eqref{eq:tail-sector-pair}
by the positive gcd of its twelve entries.

At $\rho=\one$,
\[
 I_0(\one)=\int_{[-1,1]^2}
       (3q_{0,{\rm sq}}^\alpha+9q_{0,{\rm rh}}^\alpha)\,dx\,dy,\qquad
 I_1(\one)=64I_0(\one),\qquad P_\alpha^{\rm sym}(\one)=0 .
\]
Vertex invariance and homogeneity give the first-order vanishing
in each sector, since for every $e\in\mathcal E$,
\[
 \partial_{\rho_e}P_\alpha^{\rm sym}(\one)
 =\tfrac16\sum_{f\in\mathcal E}\partial_{\rho_f}P_\alpha^{\rm sym}(\one)
 =\tfrac{10}{6}P_\alpha^{\rm sym}(\one)=0 .
\]
Put $v=(1,1,1/2,-1/4,1/2,0)$.

\begin{lemma}[Estimate of the reference moments]\label{lem:corner}
For $\alpha\ge513$,
\begin{equation}\label{eq:asym}
 \left\|\frac{\alpha(\alpha+1)}{25}F_\alpha-v\right\|_\infty<\frac1{10}.
\end{equation}
\end{lemma}

\begin{proof}
Set $s=1-\frac45q_0$ for either reference quadratic $q_0$.
By convexity and~\eqref{eq:q0},
\[
\begin{aligned}
 \operatorname*{arg\,max}_{[-1,1]^2}q_{0,{\rm sq}}
 &=\{(-1,-1),(-1,1),(1,-1),(1,1)\},\\
 \operatorname*{arg\,max}_{[-1,1]^2}q_{0,{\rm rh}}
 &=\{(-1,-1),(-1,1),(1,-1)\},\\
 \max_{[0,1]^2}\tfrac45q_{0,{\rm rh}}&=\tfrac7{10},\qquad
 \max_{\substack{(x,y)\in[-1,1]^2\\xy=0}}\tfrac45q_0=\tfrac9{10}.
\end{aligned}
\]
Thus $s\le1/20$ lies in the quadrants with a maximizing corner.
In such a quadrant write the inward distances from its corner
$c$ as $rw$ and $r(1-w)$, where $0\le w\le1$ and
$0\le r\le1/\max(w,1-w)$. Then
\[
 \begin{aligned}
 s&=\ell r-Q_\ast r^2,\qquad \ell=aw+b(1-w),\\
 \bigl(a,b,\,10Q_\ast-w^2-(1-w)^2\bigr)
 &=
 \begin{cases}
 (1/5,1/5,0),&q_0=q_{0,{\rm sq}},\\
 (1/5,1/5,-w(1-w)),&q_0=q_{0,{\rm rh}},\quad c=(-1,-1),\\
 (2/5,1/5,w(1-w)),&q_0=q_{0,{\rm rh}},\quad c=(-1,1),\\
 (1/5,2/5,w(1-w)),&q_0=q_{0,{\rm rh}},\quad c=(1,-1).
 \end{cases}
 \end{aligned}
\]
Here $\ell\ge1/5$, $0<Q_\ast\le1/10$ and
\[
 s(0)=0,\qquad s\left(\frac1{\max(w,1-w)}\right)\ge\frac1{10},
 \qquad \frac{ds}{dr}=\ell-2Q_\ast r .
\]
Thus $0<s\le1/20$ occurs only on the increasing branch, where
$dx\,dy=r\,dr\,dw$ and
\[
 \begin{aligned}
 u&=\sqrt{1-4Q_\ast s/\ell^2}\ge1/\sqrt2,\qquad
 r=\frac{2s}{\ell(1+u)}\le6s,\qquad \frac{dr}{ds}=\frac1{\ell u},\\
 0\le \frac{\ell^2}{s}r\frac{dr}{ds}-1
 &=\frac{4Q_\ast s}{\ell^2}\frac{1+2/u}{(1+u)^2}
 \le 10s(8\sqrt2-10)\le\frac{40}{3}s .
 \end{aligned}
\]
For each monomial $R$ in $F_\alpha$ one has $|R|\le1$ and
$|R(x,y)-R(c)|\le2r$. Hence
\[
 \left|\frac{\ell^2}{s}r\frac{dr}{ds}R(x,y)-R(c)\right|
 \le\frac{40}{3}s+2r\le\frac{76}{3}s .
\]
The leading integral at each corner, extended to $0\le s\le1$, is
\[
 \frac{R(c)}4
 \int_0^1\frac{dw}{(aw+b(1-w))^2}\int_0^1s(1-s)^{\alpha-1}\,ds
 =\frac{R(c)}{4ab\alpha(\alpha+1)} .
\]
The sums of $R(c)/(4ab)$ give the square and rhombus entries of
$25v$, respectively,
\[
 \frac{25}4\sum_{c\in\{\pm1\}^2}(1,c_1^2)=(25,25),\qquad
 \frac{25}4(1,-1,1,1)+\frac{25}8(2,0,2,-2)
       =25\left(\tfrac12,-\tfrac14,\tfrac12,0\right).
\]
For each quadratic $\sum_c1/(4ab)\le25$. The extended leading
integral has tail
\[
 \alpha(\alpha+1)\int_{1/20}^1s(1-s)^{\alpha-1}\,ds
  =(19/20)^\alpha(1+\alpha/20).
\]
The error on $s\le1/20$ and the two tails therefore satisfy
\[
 \begin{aligned}
 \left\|\frac{\alpha(\alpha+1)}{25}F_\alpha-v\right\|_\infty
 &\le \frac{76}{3}\alpha(\alpha+1)
               \int_0^1s^2(1-s)^{\alpha-1}\,ds\\
 &\quad+\alpha(\alpha+1)\int_{1/20}^1s(1-s)^{\alpha-1}\,ds
        +\frac{\alpha(\alpha+1)}{25}(19/20)^{\alpha-1}\\
 &\le\frac{152}{3(\alpha+2)}
       +\frac{(\alpha+1)(\alpha+25)}{25}(19/20)^{\alpha-1}.
 \end{aligned}
\]
Both terms decrease for $\alpha\ge513$, since
\[
 \frac{d}{d\alpha}\log\!\left((\alpha+1)(\alpha+25)
                                (19/20)^{\alpha-1}\right)
 =\frac1{\alpha+1}+\frac1{\alpha+25}-\log(20/19)
 <\frac2{513}-\frac1{20}<0 .
\]
As $(19/20)^4<5/6$ and $(5/6)^8<1/4$, at $\alpha=513$ the bound is
less than
\[
 \frac{152}{1545}+11062\,2^{-32}<\frac1{10}.
 \qedhere
\]
\end{proof}

\begin{proof}[Proof of Proposition~\ref{prop:large}]
For $\alpha\ge513$ and $(C,B)\in\mathcal C$,
Numerical Lemma~\ref{num:tail-forms} gives
\[
 4B\cdot v=4B_1+4B_2+2B_3-B_4+2B_5\ge\tfrac{13}{25}\|B\|_1,
 \qquad \|C\|_1\le\|B\|_1.
\]
Together with Lemma~\ref{lem:corner}, these inequalities give
\[
\begin{split}
 \frac{\alpha(\alpha+1)}{25}(C+\alpha B)\cdot F_\alpha
 &\ge (C+\alpha B)\cdot v
      -\tfrac1{10}(\|C\|_1+\alpha\|B\|_1)\\
 &\ge\left(\frac{3\alpha}{100}-\frac{11}{10}\right)\|B\|_1>0.
\end{split}
\]
Numerical Lemma~\ref{num:tail-finite} gives the same positivity
for $35/2\le\alpha\le514$. Thus every coefficient in
$P_\alpha^{\rm sym}$ after~\eqref{eq:sector} is nonnegative for
$\alpha\in[35/2,514]\cup[513,\infty)$. Since all sector variables
are nonnegative, $P_\alpha^{(1)}=P_\alpha^{\rm sym}/24\ge0$.
Lemma~\ref{lem:tail-reduction} gives the desired moment inequality
with strictness for unequal conorms. The equality case is
\eqref{eq:conorm-equality}.
\end{proof}

\begin{proof}[Completion of the proof of Theorem~\ref{thm:main}]
Propositions~\ref{prop:bounded} and~\ref{prop:large} give
\eqref{eq:minimizers} for $1\le p\le35$ and $p\ge35$,
respectively.
\end{proof}

\section{Appendix: Numerical Computations}\label{app:certificate}

The following numerical lemmas use finite rational calculations.
For intervals $[l,u]$, $[l',u']$ and a rational number $r$, use
\[
 \begin{aligned}
 [l,u]+[l',u']&=[l+l',u+u'],\\
 r[l,u]&=[\min(rl,ru),\max(rl,ru)],\\
 R_b[l,u]&=2^{-b}[\lfloor2^b l\rfloor,\lceil2^b u\rceil].
 \end{aligned}
\]
For the product of an exact polynomial and an interval polynomial,
the coefficient at multi-index $\nu$ is
\[
 \sum_{0\le\mu\le\nu} a_\mu [l_{\nu-\mu},u_{\nu-\mu}],
\]
where $a_\mu$ are the exact coefficients. This rule includes the
exponent variable when present. Missing coefficients are zero.
In~\eqref{eq:sector-coeff} one can apply
$X_i\mapsto X_i+X_{i+1}$ for $i=0,1,2,3,4$ in that order,
starting from $\rho_{\sigma(i)}=X_i$.

The two Bernstein calculations below use the identity
\[
 t^i=\sum_{j=i}^m
      \frac{\binom ji}{\binom mi}\binom mj t^j(1-t)^{m-j},
 \qquad 0\le i\le m.
\]
It follows by expanding $t^i(t+(1-t))^{m-i}$. On $[0,1]$
the degree-$m$ Bernstein basis is nonnegative and sums to one.

\subsection{Bounded exponents}
For a comparison polynomial of degree $d$, count the coefficients
of $X_0^{\nu_0}\cdots X_5^{\nu_5}$ separately in each of the
thirty orders $\rho_{\sigma(i)}=X_i+\cdots+X_5$. The total is
$30\binom{d+5}{5}$, of which $180$ have $\nu_5\ge d-1$ and
vanish by Section~\ref{sec:bounded}.

\begin{numlemma}[Quartic Taylor comparison]\label{num:taylor-low}
In $P_\alpha^{\rm T}$ with $K=4$ and $d=22$, exactly
$2\,301\,600$ coefficients exceed $5/2$ for every
$\alpha\in[1/2,1]$. The other $120\,300$ coefficients vanish
identically on this interval.
\end{numlemma}

\begin{numlemma}[Cubic Taylor comparison]\label{num:taylor-high}
In $P_\alpha^{\rm T}$ with $K=3$ and $d=18$, exactly
$960\,780$ coefficients exceed $7$ for every $\alpha\in[1,2]$.
The other $48\,690$ coefficients vanish identically on this interval.
\end{numlemma}

\begin{numlemma}[Weighted quadratic comparison]\label{num:quadratic}
For each $I\in\{[2,5/2],[5/2,3]\}$, exactly $332\,730$ coefficients
of $P_\alpha^{(2)}$ of degree $d=14$ exceed $3/4$ for every
$\alpha\in I$. The other $16\,110$ coefficients vanish identically on $I$.
\end{numlemma}

\begin{numlemma}[Weighted cubic comparison]\label{num:cubic}
For each $I\in\{[j,j+1]:j=3,\ldots,15\}\cup\{[16,35/2]\}$,
exactly $960\,780$ coefficients of $P_\alpha^{(3)}$ of degree
$d=18$ exceed $260$ for every $\alpha\in I$. The other
$48\,690$ coefficients vanish identically on $I$.
\end{numlemma}

\begin{proof}[Proof of Numerical Lemmas~\ref{num:taylor-low}--\ref{num:cubic}]
For $I=[\alpha_-,\alpha_+]$, write
\[
 \alpha=c+rw,\qquad
 c=\frac{\alpha_-+\alpha_+}{2},\qquad
 r=\frac{\alpha_+-\alpha_-}{2},\qquad |w|\le1 .
\]
For $q\in\{q_{0,{\rm sq}},q_{0,{\rm rh}}\}$, the moment is
\[
 \frac14\int_{[-1,1]^2} x^a y^b q^{c+rw-k}\,dx\,dy,
 \qquad k,a,b\in\Z_{\ge0},\quad a+b\le2k.
\]
The coefficients in
$\binom{c-k+rw}{\ell}=\sum_{j=0}^\ell b_{\ell,j}w^j$ satisfy
\[
 b_{0,0}=1,\qquad
 b_{\ell+1,j}=\frac{(c-k-\ell)b_{\ell,j}+r b_{\ell,j-1}}{\ell+1},
 \qquad b_{\ell,j}=0\quad(j<0\text{ or }j>\ell).
\]
The coefficient of $w^j$ in the moment truncated at degree $J$
in $q-1$ is
\begin{equation}\label{eq:bounded-moment-coefficient}
 \frac14\sum_{\ell=j}^{J} b_{\ell,j}
   \int_{[-1,1]^2}x^ay^b(q-1)^\ell\,dx\,dy.
\end{equation}
All integrals here are given by~\eqref{eq:multinomial-moment}.

For the Taylor comparisons, retain $0\le j\le h$ and use
\[
 (J,h)=
 \begin{cases}
 (40,10),&K=4,\\
 (50,12),&K=3.
 \end{cases}
\]
The coefficient $\ell^1$ norm of the spatial term of order
$\ell$ is bounded by
\[
 4^{-\ell}\prod_{i=0}^{\ell-1}
     \frac{|c-k-i|+r}{i+1}
 \le 4^{-\ell}\binom{\ell+3}{3},\qquad 0\le k\le4.
\]
The empty product is one. The spatial tail in the retained
powers of $w$ satisfies
\[
 \sum_{\ell=J+1}^{\infty}4^{-\ell}\binom{\ell+3}{3}
 \le
 \frac{4^{-(J+1)}\binom{J+4}{3}}{1-(J+5)/(4(J+2))}.
\]
For the omitted powers of $w$ use the exponential series for
$q^{c-k}\exp(rw\log q)$. Since $q^{c-k}\le4$ and
$r|\log q|<r/3$, the total error is bounded by
\[
 \frac{4^{-(J+1)}\binom{J+4}{3}}{1-(J+5)/(4(J+2))}
 +\frac{4(r/3)^{h+1}}{(h+1)!\,(1-r/(3(h+2)))}
 <2^{-60}.
\]

For Numerical Lemma~\ref{num:cubic}, retain $0\le j\le19$
on each stated interval. Scale the moments by $2^{-s}$, with
\[
 s=\min\{j\in\Z_{\ge1}:2^j\ge (5/4)^{\lceil c\rceil}\}.
\]
Indeed, for $0\le k\le4$,
\[
 q^{c-k}\le
 \begin{cases}
 (5/4)^{\lceil c\rceil},&c\ge k,\\
 4/3,&c<k ,
 \end{cases}
 \qquad
 4/3<(5/4)^{\lceil c\rceil}\quad(c\ge7/2).
\]
From $2^{-s}q^{c-k}\le1$ and $|\log q|<3/10$, the
normalized exponential tail is at most
\[
 \frac{(3r/10)^{20}}{20!\,(1-(3r/10)/21)}<2^{-61}.
\]
Take the least $J\ge\lceil c\rceil$ for which
\[
 \frac43\,2^{-s}\max_{0\le k\le4,\ k\in\Z}
 \left(4^{-(J+1)}
 \prod_{i=0}^{J}\frac{|c-k-i|+r}{i+1}\right)<2^{-61}.
\]
Since $c-r\ge3$, $k\le4$, and $\ell\ge\lceil c\rceil$,
\[
 \frac{|c-k-\ell|+r}{4(\ell+1)}
 =\frac{\ell+k-(c-r)}{4(\ell+1)}\le\frac14 .
\]
The factor $4/3$ sums this geometric tail. The total error is
less than $2^{-60}$.

For Numerical Lemma~\ref{num:quadratic}, use $\alpha=n+w/2$
with $0\le k\le2$ and
\[
 (n,w)\in
 \begin{cases}
 \{2\}\times[0,1],&I=[2,5/2],\\
 \{3\}\times[-1,0],&I=[5/2,3].
 \end{cases}
\]
Use $b_{\ell,j}$ from the same recurrence with $c-k=0$ and
$r=1/2$. For the moments scaled by $(4/5)^n$, the retained
coefficients are
\[
 \frac{(4/5)^n}{4}
 \sum_{\ell=j}^{36} b_{\ell,j}
 \sum_{i=0}^{n-k}\binom{n-k}{i}
     \int_{[-1,1]^2}x^a y^b(q-1)^{\ell+i}\,dx\,dy,
 \qquad 0\le j\le12 .
\]
The coefficient norm of $\binom{w/2}{\ell}$ is at most one.
Using $(4/5)^nq^{n-k}\le(4/3)^k$ and $|\log q|/2<1/6$ gives
the following bound for the two tails,
\[
 \frac{16}{9}\left(\frac{4}{3\cdot4^{37}}
                 +\frac{2}{6^{13}13!}\right)<2^{-60}.
\]
On $I=[5/2,3]$, replace $w$ by $-w$ in the power expansion
before testing the coefficients on $[0,1]$.

Parity and the substitution $y\mapsto-y$ give
\[
\begin{aligned}
 \int_{[-1,1]^2}x^ay^b q_{0,{\rm sq}}^{\alpha-k}\,dx\,dy&=0
                  \quad\text{if $a$ or $b$ is odd},\\
 \int_{[-1,1]^2}x^ay^b q_{0,{\rm rh}}(x,-y)^{\alpha-k}\,dx\,dy
 &=(-1)^b\int_{[-1,1]^2}x^ay^b q_{0,{\rm rh}}(x,y)^{\alpha-k}\,dx\,dy .
\end{aligned}
\]
If the retained coefficients of the scaled moment are $a_j$, its
interval coefficient at degree $j$ is
\[
 2^{-60}
 \bigl[\lfloor2^{60}a_j\rfloor-\one_{\{j=0\}},\
       \lceil2^{60}a_j\rceil+\one_{\{j=0\}}\bigr].
\]
Use $[0,0]$ at every degree for an odd square moment. On a
reflected rhombus apply $(-1)^b$ to the geometric weight of $x^ay^b$.

For the Taylor comparisons compute the geometric coefficients
$d_{\nu,a,b}(T,k)$ in~\eqref{eq:ik-coeff} by multiplication
of $W_T(4q_T\det Q)^k$. Insert the moment intervals separately
for every tile and sum the contributions to $I_k$. For the
reference term below use the constant-weight moments at $k=0$,
one per tile, with factor $4$ on the three square tiles and
$3$ on the twelve rhombus tiles. With the $1/4$ in the
moments, this gives the reference integral in~\eqref{eq:lowpoly}.
Evaluate the following product in the displayed order and apply
$R_{32}$ to its coefficients before~\eqref{eq:sector-coeff}.
\[
\begin{aligned}
 9216P_\alpha^{\rm T}
 ={}&
 \left[24\left(1+\frac{\alpha}{3}\right)
       \Bigl(\sum_{e\in\mathcal E}\rho_e\Bigr)^3
       +108(3-\alpha)\det Q\right]\\
 &\quad\times
 \sum_{k=0}^K\left[384\cdot3^kc_k(\alpha)
       \Bigl(2\det Q\sum_{e\in\mathcal E}\rho_e\Bigr)^{K-k}\right] I_k\\
 &-\left[24\left(1-\frac{\alpha}{3}\right)
       \Bigl(\sum_{e\in\mathcal E}\rho_e\Bigr)^3
       +108(3+\alpha)\det Q\right]\\
 &\quad\times
 \left[24\det Q\Bigl(2\det Q\sum_{e\in\mathcal E}\rho_e\Bigr)^K
   \int_{[-1,1]^2}\bigl(3q_{0,{\rm sq}}^\alpha
                          +9q_{0,{\rm rh}}^\alpha\bigr)\,dx\,dy\right].
\end{aligned}
\]

For the weighted comparisons, let $\epsilon_T=-1$ on the
reflected rhombus tiles and $\epsilon_T=1$ on the other tiles.
In~\eqref{eq:ik-coeff} the exact weight of
$\frac14\int_{[-1,1]^2}x^ay^bq^{\alpha-k}\,dx\,dy$ is
\[
 \sum_{\substack{T\in\mathcal T\\q_{0,T}(x,\epsilon_T y)=q(x,y)}}
    \epsilon_T^b
    \sum_{\substack{\nu\in\Z_{\ge0}^6\\
                    \sum_{e\in\mathcal E}\nu_e=4k+3}}
       d_{\nu,a,b}(T,k)\prod_{e\in\mathcal E}\rho_e^{\nu_e},
 \qquad q\in\{q_{0,{\rm sq}},q_{0,{\rm rh}}\}.
\]
The equality in the outer index is an identity in $x,y$.
Omit the zero square moments. The reference integral uses the two
$(k,a,b)=(0,0,0)$ moments with weights $12$ and $36$.
Multiply these weights by the exact conorm factors of
\eqref{eq:quadpoly} or~\eqref{eq:cubicpoly} before inserting
the moment intervals. The three exponent factors are
\[
 \begin{array}{ll}
 2(\alpha-1),\quad 2(\alpha-2),\quad 2
   &\text{in the quadratic comparison},\\
 4(\alpha-2),\quad 4(\alpha-3),\quad 4
   &\text{in the cubic comparison}.
 \end{array}
\]
For the quadratic comparison substitute~\eqref{eq:sector-coeff}
in the exact conorm weights and apply $R_{52}$ to the moment
coefficients before taking their products.

For the cubic comparison take the products and sum in the conorms.
At each conorm coefficient write the resulting exponent coefficients
as $[l_j,u_j]$. Retain degrees up to
\[
 h=\begin{cases}
 8,&I=[m,m+1],\quad m=3,\ldots,15,\\
 11,&I=[16,35/2].
 \end{cases}
\]
Add the interval
$\bigl[-\sum_{j>h}\max\{|l_j|,|u_j|\},\
         \sum_{j>h}\max\{|l_j|,|u_j|\}\bigr]$
to $[l_0,u_0]$. Apply $R_{32}$ to the retained coefficients before
\eqref{eq:sector-coeff}.
The total positive factors, including the moment scales, are
\[
 \begin{cases}
 9216,&P_\alpha^{\rm T},\quad K=3\ \text{or }4,\\
 2(4/5)^n,&P_\alpha^{(2)},\\
 4\cdot2^{-s},&P_\alpha^{(3)} .
 \end{cases}
\]

Recall that $\nu$ gives the powers of $X_0,\ldots,X_5$ after
$\rho_{\sigma(i)}=X_i+\cdots+X_5$. For $\nu_5\le d-2$, the
scaled coefficient lies in $\sum_{j\ge0}w^j[l_j,u_j]$.
The Taylor comparison has degree at most $14$ for $K=4$
and at most $16$ for $K=3$. Thus for $|w|\le1$ the scaled
coefficient is bounded below by
\[
 \sum_{j\ge0}\min\{w^j l_j,w^j u_j\}
 \ge l_0-\sum_{j\ge1}\max\{|l_j|,|u_j|\}.
\]
In each Taylor or weighted cubic position the rational arithmetic
gives either $l_j=u_j=0$ for every $j$ or
\begin{equation}\label{eq:bounded-sign-test}
 l_0-\sum_{j\ge1}\max\{|l_j|,|u_j|\}>
 \begin{cases}
 23\,040,&P_\alpha^{\rm T},\ K=4,\\
 64\,512,&P_\alpha^{\rm T},\ K=3,\\
 260,&P_\alpha^{(3)}.
 \end{cases}
\end{equation}
In the quadratic comparison, if all exact moment weights are zero,
the coefficient vanishes. In every other position, on the
parameter interval $0\le w\le1$, its lower bound is
\[
 \max\left\{
  l_0+\sum_{j=1}^{13}\min(0,l_j),\
  \min_{0\le i\le13}
  \sum_{j=0}^{i}\frac{\binom{i}{j}}{\binom{13}{j}}\,
       2^{-28}\big\lfloor2^{28}l_j\big\rfloor
 \right\}>1 .
\]
The first bound uses $0\le w^j\le1$. The second uses the
degree-$13$ Bernstein coefficients of the polynomial with power
coefficients $2^{-28}\lfloor2^{28}l_j\rfloor$. The exact scale
conversions are
\[
 \frac{23\,040}{9216}=\frac52,\qquad
 \frac{64\,512}{9216}=7,\qquad
 \frac1{2(4/5)^n}\ge\frac{25}{32}>\frac34,\qquad
 \frac{260}{4\cdot2^{-s}}\ge260 ,
\]
where $s\ge2$ by $(5/4)^{\lceil c\rceil}\ge(5/4)^4>2$.

For each interval, enumerate the thirty orders
from~\eqref{eq:sector} and
\[
 0\le\nu_i\le d-\sum_{j=0}^{i-1}\nu_j\quad(0\le i<5),
 \qquad \nu_5=d-\sum_{j=0}^4\nu_j.
\]
Among the positions with $\nu_5\le d-2$, the zero tests give
$120\,120$, $48\,510$, $15\,930$, and $48\,510$ zeros in the four
comparisons, respectively, on each indicated interval.
Adding the $180$ positions with $\nu_5\ge d-1$ gives the stated
zero counts. The other positions pass the strict lower-bound tests.
\end{proof}

\paragraph{An explicit calculation of one coefficient.}
Consider $P_\alpha^{(2)}$ in the sector
$\rho_{01}\ge\rho_{02}\ge\rho_{03}\ge\rho_{12}\ge\rho_{13}\ge\rho_{23}$.
Since the polynomial is homogeneous of degree $14$, setting
$X_2=1$ and $X_i=0$ for $i\ne2$ makes $X_2^{14}=1$ and every
other monomial zero. At this point, \eqref{eq:sector} gives
\[
 (\rho_{01},\rho_{02},\rho_{03},\rho_{12},\rho_{13},\rho_{23})
 =(1,1,1,0,0,0),\qquad Q=\mathrm{Id}_3.
\]
The three tiles with free-edge pairs $(01,02)$, $(01,03)$ and
$(02,03)$ have $W_T=1$ and $4q_T=1+x^2+y^2$. All other
weights vanish. The reflected hexagon parametrizations
in~\eqref{eq:ik-explicit} give
\[
 I_k=\frac34\int_{[-1,1]^2}
       (1+x^2+y^2)^k q_{0,{\rm rh}}^{\alpha-k}\,dx\,dy .
\]
Substitution in~\eqref{eq:quadpoly} gives the coefficient as
\begin{equation}\label{eq:worked-coefficient}
\begin{aligned}
 \frac14\int_{[-1,1]^2}\Bigl[
 &69\,984(\alpha-1)(1+x^2+y^2)^2 q_{0,{\rm rh}}^{\alpha-2}\\
 &-116\,640(\alpha-2)(1+x^2+y^2) q_{0,{\rm rh}}^{\alpha-1}\\
 &-46\,656q_{0,{\rm sq}}^\alpha
 -139\,968q_{0,{\rm rh}}^\alpha
 \Bigr]\,dx\,dy .
\end{aligned}
\end{equation}
For $\alpha=2+w/2$, $0\le w\le1$, use the degree-$36$
binomial and degree-$12$ power truncations for the moments
scaled by $(4/5)^2$, as in Numerical Lemma~\ref{num:quadratic}.
Round the moment intervals outward by $R_{52}$ and multiply
by twice the corresponding factors in~\eqref{eq:worked-coefficient}.
This encloses $32/25$ times the coefficient in
\eqref{eq:worked-coefficient}. The lower endpoints $l_j$ satisfy
\[
 2^{52}l_0=122\,896\,820\,705\,070\,320\,064,\qquad
 2^{52}\sum_{j=1}^{13}\min(0,l_j)
     =-962\,510\,680\,643\,051\,232 .
\]
Thus on $2\le\alpha\le5/2$ the exact coefficient is bounded below by
\[
 \frac{25}{32}\Bigl(l_0+\sum_{j=1}^{13}\min(0,l_j)\Bigr)
 =\frac{25\cdot121\,934\,310\,024\,427\,268\,832}{2^{57}}>21\,152 .
\]
At $\alpha=2$, direct integration in~\eqref{eq:worked-coefficient}
gives $106\,596/5$.

\subsection{Large exponents}

\begin{numlemma}[Degree-ten coefficient pairs]\label{num:tail-forms}
For a conorm tuple $\rho$, set
\begin{equation}\label{eq:tail-vector-polynomials}
\begin{aligned}
 C(\rho)={}&24(\det Q)\Bigl(\sum_{e\in\mathcal E}\rho_e\Bigr)
       \left[2\Bigl(\sum_{e\in\mathcal E}\rho_e\Bigr)^3+54\det Q\right]\\
 &\quad{}\times\left[16(W_{01,23}+W_{02,13}+W_{12,03})-12\det Q\right]
                      (16,4,-14,4,-4,2),\\[2pt]
 B(\rho)={}&9216\left[4\Bigl(\sum_{e\in\mathcal E}\rho_e\Bigr)^3+27\det Q\right]
         \sum_{T\in\mathcal T}W_T\Phi_T(4q_T\det Q)\\
 &-384(\det Q)\Bigl(\sum_{e\in\mathcal E}\rho_e\Bigr)
        \left[2\Bigl(\sum_{e\in\mathcal E}\rho_e\Bigr)^3+54\det Q\right]\\
 &\quad{}\times\bigl[(W_{01,23}+W_{02,13}+W_{12,03})(16,4,-14,4,-4,2)\\
 &\hspace{44mm}{}+3(\det Q)(0,0,14,-4,4,-2)\bigr].
\end{aligned}
\end{equation}
For the tile with free edges $u,v$ and normal $\eta$,
$W_T=\rho_u\rho_v\eta^TQ\eta$ and $\Phi_T$ is given by
\eqref{eq:tail-weight-vectors}. At every evaluation use
\[
 4q_T\det Q
 =\sum_{\substack{e,f\in\mathcal E\\e<f}}
         \rho_e\rho_f\det(g_e,g_f,2z_T)^2 .
\]
This polynomial expression also applies when $\det Q=0$.

For each of the thirty orders $\sigma$ and
$\nu\in\Z_{\ge0}^6$ with $\nu_0+\cdots+\nu_5=10$, set
$\rho_{\sigma(i)}(m)=m_i+\cdots+m_5$ in each summand below. Then
\begin{equation}\label{eq:tail-sector-pair}
\begin{aligned}
 (C_{\sigma,\nu},B_{\sigma,\nu})
 ={}&\frac{1}{\nu_0!\cdots\nu_5!}
 \sum_{m_0=0}^{\nu_0}\cdots\sum_{m_5=0}^{\nu_5}
 (-1)^{10-m_0-\cdots-m_5}\\
 &\qquad{}\times\left(\prod_{i=0}^5\binom{\nu_i}{m_i}\right)
                \bigl(C(\rho(m)),B(\rho(m))\bigr).
\end{aligned}
\end{equation}
There are $86\,172$ nonzero pairs, $3\,738$ zero pairs with
$\nu_5\le8$, and $180$ zero pairs with $\nu_5\ge9$.

Divide each nonzero pair by the positive gcd of its twelve
entries and retain each resulting pair $(C,B)$ once.
There are $30\,812$ such pairs, and each satisfies
\[
 \begin{gathered}
 0<\sum_{\ell=1}^6|B_\ell|,\qquad
 \sum_{\ell=1}^6|C_\ell|\le\sum_{\ell=1}^6|B_\ell|,\\
 4B_1+4B_2+2B_3-B_4+2B_5\ge\frac{13}{25}\sum_{\ell=1}^6|B_\ell|.
 \end{gathered}
\]
\end{numlemma}

\begin{numlemma}[Finite large-exponent comparison]\label{num:tail-finite}
For every pair $(C,B)$ obtained by dividing a nonzero pair
in~\eqref{eq:tail-sector-pair} by its positive gcd, and every
$35/2\le\alpha\le514$,
\[
\begin{aligned}
 \frac{(4/5)^{\alpha-1}}{4}\biggl\{&
 \int_{[-1,1]^2}
   \bigl[C_1+\alpha B_1+(C_2+\alpha B_2)x^2\bigr]
                           q_{0,{\rm sq}}^{\alpha-1}\,dx\,dy\\
 &+\int_{[-1,1]^2}
   \bigl[C_3+\alpha B_3+(C_4+\alpha B_4)x
     +(C_5+\alpha B_5)x^2\\
 &\hspace{42mm}{}+(C_6+\alpha B_6)xy\bigr]
                          q_{0,{\rm rh}}^{\alpha-1}\,dx\,dy\biggr\}
 >\frac{47}{100}.
\end{aligned}
\]
\end{numlemma}

\begin{proof}[Proof of Numerical Lemma~\ref{num:tail-forms}]
The identity
\[
 16U(\rho)-(\det Q)U(\one)
 =\left[16(W_{01,23}+W_{02,13}+W_{12,03})-12\det Q\right]
     (16,4,-14,4,-4,2)
\]
gives $P_\alpha^{\rm sym}=(C(\rho)+\alpha B(\rho))\cdot F_\alpha$
in~\eqref{eq:tail-sum}. Each polynomial in~\eqref{eq:tail-vector-polynomials}
is homogeneous of degree ten. For nonnegative integers $n,k$,
\[
 \sum_{m=0}^n(-1)^{n-m}\binom{n}{m}m^k
 =\begin{cases}
   0,& k<n,\\
   n!,& k=n.
  \end{cases}
 \qquad (k\le n,\ 0^0=1).
\]
Indeed, the left side is the $n$th forward difference of $x^k$ at
zero, and each difference lowers the degree by one and multiplies
the leading coefficient by the previous degree. A monomial
$X_0^{k_0}\cdots X_5^{k_5}$ of degree ten with $(k_0,\ldots,k_5)\ne\nu$
has $k_i<\nu_i$ for some $i$, so it contributes zero
to~\eqref{eq:tail-sector-pair}. The monomial with exponents $\nu$
contributes its coefficient times $\nu_0!\cdots\nu_5!$ before division.
This proves the formula for the two coefficient vectors.

For the numerical counts, in each order evaluate
\eqref{eq:tail-vector-polynomials} at $\rho_{\sigma(i)}=m_i+\cdots+m_5$
for the $\binom{16}{6}=8\,008$ tuples of nonnegative integers
with $m_0+\cdots+m_5\le10$. Apply~\eqref{eq:tail-sector-pair} to
the $30\binom{15}{5}=90\,090$ positions. The exact zero tests give
\[
 90\,090=86\,172+3\,738+180 .
\]
The gcd division gives $30\,812$ distinct pairs.
For each of them the stated inequalities follow by summing
the six absolute values in each vector and testing the integer
$100B_1+100B_2+50B_3-25B_4+50B_5-13\sum_{\ell=1}^6|B_\ell|$,
which is nonnegative.
\end{proof}

\paragraph{An explicit calculation of the coefficient pair for $X_2^{10}$.}
Use $\sigma=(01,02,03,12,13,23)$ and $\nu=(0,0,10,0,0,0)$. The coefficient vectors of $X_2^{10}$
are the values of the homogeneous degree-ten polynomials
$C(\rho)$ and $B(\rho)$ at $X_2=1$ and $X_i=0$ for $i\ne2$.
By~\eqref{eq:sector}, this gives
\[
 (\rho_{01},\rho_{02},\rho_{03},\rho_{12},\rho_{13},\rho_{23})
 =(1,1,1,0,0,0),\qquad Q=\mathrm{Id}_3 .
\]
Only the three rhombus tiles of weight one remain, each with
$4q_T=1+x^2+y^2$. Hence
\[
 U(\rho)=(0,0,42,-12,12,-6),\qquad
 \sum_{T\in\mathcal T}W_T\Phi_T(4q_T\det Q)=(0,0,3,0,6,0).
\]
Inserting these values and $\sum_{e\in\mathcal E}\rho_e=3$
in~\eqref{eq:tail-vector-polynomials} gives
\begin{equation}\label{eq:tail-worked-pair}
 (C_{\sigma,\nu},B_{\sigma,\nu})
 =186\,624\bigl((-8,-2,7,-2,2,-1),\ (0,0,-8,8,32,4)\bigr).
\end{equation}
The positive gcd of the twelve entries is $186\,624$.
Dividing both vectors by this number gives the pair in $\mathcal C$.

\begin{proof}[Proof of Numerical Lemma~\ref{num:tail-finite}]
Use the intervals
\[
\begin{aligned}
 \mathcal I={}&\{[35/2,18]\}\cup
       \{[j,j+1]:j=18,\ldots,57\}\\
      &{}\cup\{[58+8j,66+8j]:j=0,\ldots,56\},
\end{aligned}
\]
whose union is $[35/2,514]$. Fix
$I=[a,b]\in\mathcal I$ and put $n=\lfloor a\rfloor-1$. Then
\[
 \alpha-n-1\in
 \begin{cases}
 [1/2,1],&I=[35/2,18],\\
 [0,1],&I=[j,j+1],\quad j=18,\ldots,57,\\
 [0,8],&I=[58+8j,66+8j],\quad j=0,\ldots,56 .
 \end{cases}
\]
Recall that $(F_\alpha)_\ell=\frac14\int_{[-1,1]^2}R u^{\alpha-1}\,dx\,dy$,
with
\[
 u=\begin{cases}
 (4/5)q_{0,{\rm sq}},&\ell=1,2,\\
 (4/5)q_{0,{\rm rh}},&\ell=3,4,5,6,
 \end{cases}
 \qquad
 R=1,x^2,1,x,x^2,xy\quad\text{in order of }\ell .
\]
Since $|R|\le1$ and $0\le1-u\le2/5$, the binomial series gives
\[
 \left|(F_\alpha)_\ell-
 \frac14\sum_{j=0}^{30}(-1)^j\binom{\alpha-n-1}{j}
       \int_{[-1,1]^2}R u^n(1-u)^j\,dx\,dy\right|
 \le\varepsilon_\ell ,
\]
where
\[
 \varepsilon_\ell=\frac{(2/5)^{31}}{1-2/5}
 \begin{cases}
 1,& I\subseteq[35/2,58],\\[2pt]
 \dfrac{8^8}{4\cdot8!}\displaystyle\int_{[-1,1]^2}u^n\,dx\,dy,
   & I\subseteq[58,514].
 \end{cases}
\]
For the first case, the product for $\binom{\alpha-n-1}{j}$
has absolute value at most one. In the second case,
\[
 \begin{gathered}
  \left|\binom{\alpha-n-1}{8}\right|\le\frac{8^8}{8!},\\
  \left|\binom{\alpha-n-1}{j+1}\right|
  =\frac{|j-(\alpha-n-1)|}{j+1}
    \left|\binom{\alpha-n-1}{j}\right|
  \le\left|\binom{\alpha-n-1}{j}\right|\qquad(j\ge8).
 \end{gathered}
\]
Summing the geometric tail gives the displayed bound.

For $\alpha=a+(b-a)t$, $0\le t\le1$, and $r=0,1$, the
power coefficients $a_{\ell,i}^{(r)}$ are given by
\[
 \sum_{i=0}^{31}a_{\ell,i}^{(r)}t^i
 =\frac{(a+(b-a)t)^r}{4}\sum_{j=0}^{30}(-1)^j
      \binom{a+(b-a)t-n-1}{j}
      \int_{[-1,1]^2}R u^n(1-u)^j\,dx\,dy .
\]
This polynomial differs from $\alpha^r(F_\alpha)_\ell$ by at
most $b^r\varepsilon_\ell$. The spatial integrals are
\[
 \int_{[-1,1]^2}R u^n(1-u)^j\,dx\,dy
 =\sum_{i=0}^{j}(-1)^i\binom ji\int_{[-1,1]^2}R u^{n+i}\,dx\,dy,
 \qquad 0\le j\le30.
\]
For $I\subseteq[35/2,58]$ and $R=x^ay^b$,
\eqref{eq:multinomial-moment} gives
$\frac14\int_{[-1,1]^2}Ru^{n+i}\,dx\,dy$ by using $u$ as the
quadratic and $n+i$ as its power.
For $I\subseteq[58,514]$, use the coordinates $Y$ and domains
$\Omega$ given by
\[
\begin{aligned}
 u=\tfrac45q_{0,{\rm sq}}:\quad
 &Y=(x,y),\quad \Omega=[-1,1]^2,\quad
 u=\tfrac45+\tfrac{Y_1^2+Y_2^2}{10},\\
 u=\tfrac45q_{0,{\rm rh}}:\quad
 &Y=(x-1,y-1),\quad \Omega=[-2,0]^2,\quad
 u=\tfrac35+\tfrac{Y_1^2-Y_1Y_2+Y_2^2}{10}.
\end{aligned}
\]
The translated weights are $1$, $Y_1+1$, $(Y_1+1)^2$ and
$(Y_1+1)(Y_2+1)$. Since $Y\cdot\nabla u=2(u-u(0))$, the
divergence theorem gives, for $m\ge1$,
\[
 \begin{aligned}
 (h+k+2+2m)\int_\Omega Y_1^hY_2^k u^m\,dY
 &=\int_{\partial\Omega}(Y\cdot\nu)Y_1^hY_2^k u^m\,dS\\
 &\quad{}+2m u(0)\int_\Omega Y_1^hY_2^k u^{m-1}\,dY,\\
 \int_{[L,H]^2}Y_1^hY_2^k\,dY
 &=\frac{(H^{h+1}-L^{h+1})(H^{k+1}-L^{k+1})}{(h+1)(k+1)} .
 \end{aligned}
\]
Here $\nu$ is the outward unit normal, and the last line gives
the initial values at $m=0$.
Every edge with $Y\cdot\nu\ne0$ has $u=(9+s^2)/10$ for
$-1\le s\le1$. Define
\[
 K_{j,m}=\int_{-1}^1 s^j\left(\frac{9+s^2}{10}\right)^m ds .
\]
Integrating the derivative of $s^{j+1}((9+s^2)/10)^m$ gives
\[
 \begin{aligned}
 K_{j,0}&=\frac{1-(-1)^{j+1}}{j+1},\\
 (j+1+2m)K_{j,m}&=1-(-1)^{j+1}+\frac{9m}{5}K_{j,m-1}
 \qquad(m\ge1).
 \end{aligned}
\]
On the square take $s=Y_j$ on the edge with fixed coordinate
$Y_i=\pm1$. On the translated rhombus the edges $Y_i=0$
have $Y\cdot\nu=0$, and on $Y_i=-2$ take $s=Y_j+1$.
The boundary term for $Y_1^hY_2^k$ is therefore
\[
\begin{cases}
 (1+(-1)^h)K_{k,m}+(1+(-1)^k)K_{h,m},
    &\Omega=[-1,1]^2,\\[2pt]
 2(-2)^h\displaystyle\sum_{j=0}^k\binom{k}{j}(-1)^{k-j}K_{j,m}
 +2(-2)^k\displaystyle\sum_{j=0}^h\binom{h}{j}(-1)^{h-j}K_{j,m},
    &\Omega=[-2,0]^2 .
\end{cases}
\]
For $j=0,\ldots,31$, add the residual interval
$[-b^r\varepsilon_\ell,b^r\varepsilon_\ell]$ to each Bernstein
coefficient of the approximating polynomial. The resulting
integer endpoints at scale $2^{60}$ are
\[
 \begin{aligned}
 l_{\ell,j}^{(r)}
 &=\left\lfloor2^{60}
       \left(\sum_{i=0}^j\frac{\binom ji}{\binom{31}i}
                          a_{\ell,i}^{(r)}-b^r\varepsilon_\ell\right)
       \right\rfloor,\\
 u_{\ell,j}^{(r)}
 &=\left\lceil2^{60}
       \left(\sum_{i=0}^j\frac{\binom ji}{\binom{31}i}
                          a_{\ell,i}^{(r)}+b^r\varepsilon_\ell\right)
       \right\rceil .
 \end{aligned}
\]
The Bernstein basis is nonnegative and sums to one on $[0,1]$, so
\[
 \alpha^r(F_\alpha)_\ell\in 2^{-60}\sum_{j=0}^{31}
 \binom{31}{j}t^j(1-t)^{31-j}[l_{\ell,j}^{(r)},u_{\ell,j}^{(r)}] .
\]
Using the same basis for both terms of $(C+\alpha B)\cdot F_\alpha$
gives the lower bound
\[
 \frac1{2^{60}}\min_{0\le j\le31}\sum_{\ell=1}^{6}
 \left[
   \min\{C_\ell l_{\ell,j}^{(0)},C_\ell u_{\ell,j}^{(0)}\}
  +\min\{B_\ell l_{\ell,j}^{(1)},B_\ell u_{\ell,j}^{(1)}\}
 \right].
\]
Taking the minimum of this expression over $(C,B)\in\mathcal C$
and $I\in\mathcal I$ with $I\subseteq J$ gives
\[
 \begin{cases}
 \dfrac{4\,288\,754\,204\,510\,563\,386}{2^{60}},
       &J=[35/2,58],\\[3pt]
 \dfrac{552\,536\,073\,955\,686\,419}{2^{60}},
       &J=[58,514].
 \end{cases}
\]
Both minima exceed $47/100$. They occur at $j=31$ on
$[57,58]$ and $[506,514]$, respectively, for
\[
 (C,B)=\bigl((-8,-2,7,-2,2,-1),\ (0,0,-8,8,32,4)\bigr).\qedhere
\]
\end{proof}

\noindent
{\bf Acknowledgements.}
 The research of S. Luo is partially supported by the National Natural Science Foundation of China (NSFC) under Grant Nos. 12261045 and 12001253, and by the Jiangxi Jieqing Fund under Grant No. 20242BAB23001.
 The author used AI to assist the long but finite computations.
\vskip0.1in
{\bf Statements and Declarations: there is no conflict of interest.}

{\bf Data availability: the manuscript has no associated data.}

\end{document}